\documentclass{article}
\usepackage{amsmath,amssymb,amsthm,a4wide,url}
\usepackage[onehalfspacing]{setspace}
\usepackage[numbers]{natbib}
\newcommand{\nw}{}
\newcommand{\old}{}
\newcommand{\cn}{\,{:}\;}
\newtheorem{theorem}{Theorem}[section]
\newtheorem{proposition}[theorem]{Proposition}
\newtheorem{lemma}[theorem]{Lemma}
\newtheorem{corollary}[theorem]{Corollary}
\theoremstyle{definition}
\newtheorem{definition}[theorem]{Definition}
\newtheorem{example}[theorem]{Example}
\newtheorem{remark}[theorem]{Remark}

\newcommand{\mC}{\mathbb{C}}

\newcommand{\cH}{\mathcal{H}}
\newcommand{\cM}{\mathcal{M}}
\newcommand{\bY}{\overline{Y}}

\newcommand{\hx}{\hat{x}}

\newcommand{\ty}{\tilde{y}}
\newcommand{\tH}{\tilde{H}}

\newcommand{\tM}{\tilde{M}}
\newcommand{\tT}{\tilde{T}}
\newcommand{\tV}{\tilde{V}}
\newcommand{\tY}{\tilde{Y}}
\newcommand{\tZ}{\tilde{Z}}
\newcommand{\hM}{\hat{M}}
\newcommand{\hV}{\hat{V}}
\newcommand{\hY}{\hat{Y}}
\newcommand{\oT}{\overline{T}}
\newcommand{\hl}{\hat{\ell}}
\DeclareMathOperator{\im}{im}
\DeclareMathOperator{\diag}{diag}
\begin{document}
\title{\textbf{Keldysh's theorem revisited}}
\author{Johannes M. Schumacher\thanks{%
Amsterdam School of Economics, University of Amsterdam, Amsterdam, the 
Netherlands. E-mail: \protect\url{j.m.schumacher@uva.nl}. ORCID 0000-0001-5753-3412.}}
\maketitle
\begin{abstract}
\noindent
In a variety of applications, the problem comes up of describing the principal part of the inverse of a holomorphic 
operator
at an eigenvalue in terms of left and right root functions associated to the eigenvalue. Such a description was given
by Keldysh in 1951. His theorem, the proof of which was published only in 1971, is a fundamental result in
the local spectral theory of operator-valued functions. Here we present a
streamlined derivation in the matrix case, and we extend Keldysh's theorem by means of
a new principal part formula. Special attention is given to the semisimple case (first-order poles).
\end{abstract}

\section{Introduction}

In 1951, M.V.\ Keldysh gave a representation theorem for the principal part, at a given point in the  complex
plane, of operator-valued functions of the form $(I+A(\lambda))^{-1}$, where $A(\lambda)$ is a compact linear 
operator that depends holomorphically on the complex parameter $\lambda$. 
Keldysh's work was primarily motivated by applications 
to the dynamic analysis of flexible structures, in particular aircraft. Similar applications have motivated
another early stream of research, exemplified by \cite{Frazer} and \cite{Lancaster}, which concentrated 
on matrices (rather than operators) that depend on a complex parameter.\footnote{Holomorphic
operator-valued functions are classically used to describe relations between forces and displacements in a 
mechanical structure
as a function of a complex (Fourier transform) parameter. The matrix case appears when attention is focused on forces 
and displacements at finitely many points; the vibration of
the structure itself may then still be described by a partial differential equation.} In more recent years, while the local 
spectral theory of 
matrix-valued functions has remained highly relevant in the study of vibrating structures, applications in
other areas have emerged as well. Perhaps most notably, a key result in \cite{Engle}, one of the publications for 
which the
authors received the Nobel Prize in Economics in 2003, is a representation theorem for
the principal part of the inverse of a holomorphic matrix at a particular point of interest in the complex plane. 
In the further development of the literature on cointegration in econometric theory, the local 
spectral theory of matrix-valued functions has continued to play a central role; see for example 
\cite{laCour}, \cite{Johansen2008}, and \cite{Faliva}. For the use of Keldysh's theorem in 
the numerical analysis of nonlinear eigenvalue problems, one may consult for instance \cite{MM} and \cite{Guttel}.
Kozlov and Maz'ya \cite[p.\,6]{Kozlov}, in their monograph on differential equations with operator
coefficients, refer to the theorem as a fundamental result. Generally speaking,
Keldysh's theorem can play a key role in situations where transform analysis is applied in a multivariate
setting; for recent examples in the context of Markov-modulated stochastic processes, see
\cite{Beare2024} and \cite{Beare}.

Keldysh's 1951 paper is a short note without proofs, containing the representation theorem as well as a completeness
theorem for extended systems of eigenvectors. Proofs of both theorems
were provided by Keldysh in a manuscript that was written in 1950 and circulated in 1951, but that was published only
in 1971.\footnote{See
the editor's note in \cite{Keldysh1971}.} Later authors have given alternative proofs of (versions of) Keldysh's 
representation theorem; see
for instance \cite{GS}, \cite{Mennicken}. Some partial results have been rediscovered independently, as
discussed later in this paper. Textbook treatments can be found in \cite{MM} and
\cite{Kozlov}, both in an operator setting. 

For many applications, the matrix case is sufficient. In this case,
the theory has its own characteristics; techniques from commutative algebra can be used with ease, and 
it is natural to treat left and right eigenvectors on an equal footing. We present below a derivation of 
Keldysh's theorem in the matrix case using an
approach that might be termed geometric, since the treatment is centered around a certain
sequence of subspaces associated to a singularity of a matrix-valued function. \nw{A running example
is used to illustrate all essential concepts and results.} The classical theory may in this way
become more accessible to researchers interested in matrix-valued functions. 

Broadly (and vaguely) speaking, the purpose of Keldysh's theorem and related results is to use the ``local zero structure''
of a nonsingular holomorphic matrix/operator to describe the ``local pole structure'' of its inverse.
The following special case may serve as an illustration. Consider a square matrix $T(\lambda)$ of size $n \times n$ whose 
elements are holomorphic functions of the complex variable $\lambda$, and which is such that $\det T(\lambda)$ is not identically
zero. Assume that there is a singularity at $\lambda = \lambda_0$, i.e., $\det T(\lambda_0) = 0$. The constant matrix $T(\lambda_0)$
then has nontrivial left and right nullspaces with dimension, say, $r$; so, there exist matrices $Z_L \in \mC^{r \times n}$
and $Z_R \in \mC^{n \times r}$, of full row and column rank respectively, such that $Z_L T(\lambda_0) = 0$ and
$T(\lambda_0) Z_R = 0$. Then (see Prop.\,\ref{residuethm}) the pole of the meromorphic matrix $T^{-1}(\lambda)$ at $\lambda_0$ is
of first order if and only if the constant matrix $Z_R T'(\lambda_0) Z_L$ is nonsingular, and in that case the residue of $T^{-1}(\lambda)$
at $\lambda_0$ is given by $\text{\rm Res}(T^{-1};\lambda_0) = Z_R(Z_LT'(\lambda_0)Z_R)^{-1}Z_L$. The main result of
the present paper (Thm.\,\ref{mainthm}) gives an analogous formula for the principal part of the inverse of a holomorphic
matrix at a singularity of arbitrary order. 

\section{Problem setting}

Let $\Omega$ be a fixed open subset of the complex plane. The ring of holomorphic functions and the field
of meromorphic functions defined on $\Omega$ are denoted by $\cH$ and by $\cM$, respectively. 
Throughout the paper, attention is focused on a fixed point $\lambda_0 \in \Omega$.
The ring of locally holomorphic functions (i.e., functions that are holomorphic in a neighborhood of 
$\lambda_0$) is denoted by $\cH_0$. Due to the fact that zeros of holomorphic functions are isolated, 
a function $y \in \cH_0$ is a \emph{unit} in $\cH_0$ (i.e., has an inverse in $\cH_0$) if and only if 
$y(\lambda_0) \neq 0$. 

A matrix $M \in \cM^{n \times n}$ is said to be \emph{invertible} when it has an inverse 
in $\cM^{n \times n}$, or, equivalently, when $\det M$ is not identically zero. A matrix $M \in \cH_0^{n\times n}$ 
that has an inverse in $\cH_0^{n \times n}$ is said to be \emph{unimodular} (with respect to the ring $\cH_0$).
It follows from Cramer's rule that $M \in \cH_0^{n\times n}$ is unimodular if and only if $\det M$ is a 
unit in $\cH_0$, i.e., if and only if the constant matrix $M(\lambda_0) \in \mC^{n\times n}$ is nonsingular.
We will say that a matrix $M \in \cH_0^{n \times m}$ is \emph{left (right) unimodular} when it has a left
(right) inverse in $\cH_0^{m \times n}$; this happens if and 
only if the constant matrix $M(\lambda_0)$ has full column (row) rank.

For notational convenience, define $\chi_0 \in \cH$ by 
\begin{equation}
\chi_0(\lambda) = \lambda - \lambda_0 \qquad\qquad (\lambda \in \Omega).
\end{equation}
The Laurent expansion around $\lambda_0$ of a matrix $M \in \cM^{n \times m}$ can then be written as
\begin{equation} \label{Laurent}
M = \sum_{j = -s}^\infty M_j \chi_0^j
\end{equation}
with $M_j \in \mC^{n \times n}$ ($j=-s,-s+1,\dots$). The meromorphic matrix $M$ is said to have 
a \emph{pole} at $\lambda_0$ if there exists $j<0$ such that $M_j \neq 0$, and the \emph{order} of the pole at
$\lambda_0$ is the smallest value of $j$ for which this holds. Matrices in $\cH_0^{n \times n}$ are said to have
pole order 0 at $\lambda_0$. The coefficient 
of $\chi_0^{-1}$ in the Laurent series expansion is called the \emph{residue} of $M$ at $\lambda_0$. The part of the
sum in (\ref{Laurent}) that corresponds to indices $j < 0$ is called the
\emph{principal part} of $M$ at $\lambda_0$. The same terminology may be used when $m=1$ or $n=1$ and $M$ is
viewed as a column or row vector, rather than as a matrix. 

The notation ``\,$\doteq$\,'' is used in this paper to express \emph{equality of principal parts at $\lambda_0$}. 
In other words, given $M_1, M_2 \in \cM^{n \times m}$, we define
\begin{equation} \label{doteq}
M_1 \doteq M_2 \quad \Leftrightarrow \quad \text{there exists } H \in \cH_0^{n \times m} \text{ such that }
M_2 = M_1 + H.
\end{equation}
This could also be written in a more standard way as $M_1 = M_2 \!\mod \cH_0^{n \times m}$, but it is convenient to 
use a shorter form. The notation in (\ref{doteq}) is applied in the same way to vectors in $\cM^n$.

To avoid repetitions, throughout the paper the symbol $T$ will denote a matrix in
$\cH^{n\times n}$ that is invertible as a meromorphic matrix and that has a singularity at $\lambda_0$, i.e., 
we have $\det T(\lambda_0) = 0$ while $\det T$ is
not identically 0.  It follows from Cramer's rule that $T^{-1}$ is meromorphic.
The order of the pole of $T^{-1}$ at $\lambda_0$ is denoted by $s$, i.e., 
\begin{equation} \label{sdef}
s = \min \{ j \geq 0 \mid \chi_0^j T^{-1} \in \cH_0^{n \times n} \}.
\end{equation}
We also write
\begin{equation} \label{rdef}
r = \dim \ker T(\lambda_0).
\end{equation}
The number $r$ is called the \emph{geometric multiplicity} of the root of the holomorphic matrix $T$ at 
the point $\lambda_0$. The \emph{algebraic multiplicity} of the root of $T$ at $\lambda_0$ is the multiplicity
of $\lambda_0$ as a root of the scalar function $\det T$.

\section{Partial multiplicities} \label{Smithform}

In the ring $\cH_0$, the product of two elements can only be 0 if at least one of them is 0, and an element is divisible by
another if and only if the multiplicity of the zero at $\lambda_0$ of the first is equal to or higher than the 
multiplicity of the zero at $\lambda_0$ of the second. These properties imply in particular that $\cH_0$ is a
principal ideal domain. One can therefore use  the canonical form for matrices over
principal ideal domains that was given by H.J.S.\,Smith in 1861  (see for instance \cite[Thm.\,26.2]{MacDuffee}). 

\begin{theorem}[Smith form w.r.t.\ $\cH_0$] \label{Smith}
There exist unimodular matrices 
$U_L, U_R \in \cH_0^{n \times n}$ and uniquely determined nonnegative integers $m_1, \dots, 
m_n$ such that
\begin{equation} \label{can}
U_L T U_R = \diag(\chi_0^{m_1},\dots,\chi_0^{m_n}) \qquad\qquad (m_1 \geq \cdots \geq m_n).
\end{equation}
\end{theorem}
\vskip2mm\noindent
The numbers $m_i$ are called the \emph{partial multiplicities} of the root of $T$ at $\lambda_0$. 
It follows from (\ref{can}) that the algebraic multiplicity is the sum of the partial multiplicities. Since the pole order
of $T^{-1}$ at $\lambda_0$ is the same as that of any matrix obtained from $T^{-1}$ by left and right unimodular transformations,
we have $s = \max_{1 \leq i \leq n} m_i = m_1$. Note also that $m_r > 0$, and $m_i=0$ for $i=r+1,\dots,n$.
It will be convenient to define  a diagonal matrix containing only the diagonal elements in the Smith form that 
correspond to positive partial multiplicities:
\begin{equation} \label{deltadef}
\Delta = \diag(\chi_0^{m_1}, \dots, \chi_0^{m_r}) \in \cH_0^{r \times r}.
\end{equation} 

There is a constructive algorithm by which any given holomorphic matrix can be reduced to Smith form. The 
algorithm is laborious, however, and gives little insight in the way in which the invariants $m_i$ relate to
the coefficients in the power series expansion of $T$ around $\lambda_0$. Fortunately, an alternative
way of determining the partial multiplicities is available. Define a sequence of subspaces of $\mC^n$ as follows:%
\footnote{We will also write $L_j(T)$ in cases where $T$ needs to be specified. Mennicken and M\"{o}ller \cite[p.\,14]{MM} give a different 
but equivalent definition.}
\begin{equation} \label{LjT}
L_j = \{ y(\lambda_0) \mid y \in \cH_0^n,\; Ty  \in \chi_0^j \cH_0^n \} \qquad\qquad (j=0,1,2,\dots).
\end{equation}
Note that $L_j \supset L_{j+1}$ for all $j$, since $\chi_0^{j+1}\cH_0^n \subset \chi_0^j\cH_0^n$ for all $j$.
Also, $L_1 = \ker T(\lambda_0)$, and $L_j = \{0\}$ for $j > s$.
Define a corresponding nonincreasing sequence of nonnegative integers by
\begin{equation} \label{ljdef}
\ell_j = \dim L_j \qquad\qquad (j=0,1,2,\dots).
\end{equation}

\begin{lemma} \label{invar}
The integers $\ell_j$ defined in (\ref{ljdef}) are invariants under left and right multiplication of $T$ by 
unimodular matrices.
\end{lemma}

\begin{proof}
Fix $j \geq 0$.
Invariance under left multiplication by unimodular matrices is immediate. To prove the other half of the claim,
let a unimodular matrix $U \in \cH_0^{n\times n}$ be 
given, and take $y^0 \in L_j(TU)$. There exists $\ty \in \cH_0^n$ such that $\ty(\lambda_0) = y^0$ and 
$TU\ty \in \chi_0^j \cH_0^n$. Defining $y \in \cH_0^n$ by $y = U\ty$, we have $Ty \in \chi_0^j \cH_0^n$,
so that $y(\lambda_0) \in L_j(T)$. Since $y^0 = \ty(\lambda_0) = U^{-1}(\lambda_0)y(\lambda_0)$, this shows that
$L_j(TU) \subset U^{-1}(\lambda_0)L_j(T)$. Similarly we have $L_j(T) \subset U(\lambda_0)L_j(TU)$, and it
follows that $\dim L_j(T) = \dim L_j(TU)$. 
\end{proof}
\noindent
By reduction to Smith form, one finds the following.
\begin{proposition} \label{jl}
For all $j \geq 0$, we have
\begin{equation} \label{dimL}
\ell_j  = \#\{ 1 \leq i \leq n \mid m_i \geq j \}
\end{equation}
where the integers $m_i$ are given by (\ref{can}).
\end{proposition}
\noindent
Because the integers $m_i$ are in nonincreasing order, we can also write $\ell_j = \max\{ i \geq 1 \mid m_i \geq j \}$ 
instead of (\ref{dimL}). Therefore, $m_i \geq j$ if and only if $\ell_j \geq i$. This shows that
$(m_i)_{i=1,2,\dots}$ relates to $(\ell_j)_{j=1,2,\dots}$ in the same way as $(\ell_j)_{j=1,2,\dots}$ relates to 
$(m_i)_{i=1,2,\dots}$. Consequently, we have
\begin{equation} \label{ki}
m_i = \max\{ j \geq 1 \mid \ell_j \geq i \} = \#\{ j \geq 1 \mid \ell_j \geq i \} \qquad\qquad (i=1,\dots,n).
\end{equation}
From (\ref{dimL}) it also follows that $\#\{i \mid m_i = j \} = \ell_j - \ell_{j+1}$, so that
\begin{equation} \label{summi}
\sum_{i=1}^r m_i = \sum_{j=1}^s j (\ell_j - \ell_{j+1}) = \sum_{j=1}^s \ell_j.
\end{equation}
The sequences $(\ell_1,\dots,\ell_s)$ and $(m_1,\dots,m_r)$ are \emph{conjugate partitions} of their common sum
in the sense of enumerative combinatorics; see for instance \cite[Def.\,1.8]{Andrews}.

From (\ref{ki}), one sees that the partial multiplicities $m_i$ can be determined by computing the dimensions
of the subspaces $L_j$ ($j=1,\dots,s$).\footnote{Consequently, one might \emph{define} the partial multiplicities
$m_i$ by (\ref{ki}), rather than via the Smith form. This is the route taken in \cite{MM}. 
In the work of Keldysh \cite{Keldysh1951,Keldysh1971}, partial multiplicities appear as the maximal multiplicities that 
appear in systems of root functions whose evaluations at $\lambda_0$
form a basis for the nullspace of $T(\lambda_0)$; compare Prop.\,\ref{canchar} below.} The following proposition shows 
that finding these dimensions is a 
matter of solving finite systems of linear equations stated in terms of the first $s$ coefficients in the power series 
expansion of $T$ around $\lambda_0$.

\begin{proposition} \label{Ljchar}
For $j \geq 1$, a vector $y^0 \in \mC^n$ belongs to the subspace $L_j$ if and only if there exist vectors 
$y^1,\dots,y^{j-1} \in \mC^n$ such that
\begin{equation} \label{rooteqns}
\sum_{p=0}^k \frac{1}{p!}\,T^{(p)}(\lambda_0)y^{k-p} = 0 \quad \text{for all } 0 \leq k \leq j-1.
\end{equation}
\end{proposition}

\begin{proof}
For a vector $y \in \cH_0^n$ with power series expansion $y = \sum_{j=0}^\infty y^j \chi_0^j$ around $\lambda_0$, the
left hand side of (\ref{rooteqns}) gives the $k$-th coefficient of the power series expansion of $Ty$ around
$\lambda_0$. Therefore, if vectors 
$y^1,\dots,y^{j-1}$ are given such that (\ref{rooteqns}) is satisfied, then the function defined by 
$y = \sum_{p=0}^{j-1} \chi_0^p y^p$ satisfies $y(\lambda_0) = y^0$ and $Ty \in \chi_0^j \cH_0^n$.
Conversely, take
$y^0 \in L_j$, and let $y \in \cH_0^n$ be such that $Ty \in \chi_0^j \cH_0^n$. Since $(Ty)^{(k)}(\lambda_0) = 0$
for $0 \leq k \leq j-1$, the vectors defined 
for $p=1,\dots,j-1$ by $y^p = y^{(p)}(\lambda_0)/p!$ satisfy (\ref{rooteqns}). 
\end{proof}

\begin{example} \label{xmpl}
We illustrate the computation of the partial multiplicities via the indices $\ell_j$ in a simple example. The same example will
be used throughout the paper to illustrate various concepts and results. Take
\begin{equation} \label{Tdef}
T(\lambda) = \begin{bmatrix} \lambda^2 & \lambda^2 & 0 \\ \lambda^2 & \lambda & \lambda \\ \lambda & 0 & 1
\end{bmatrix}, \qquad \lambda_0 = 0.
\end{equation}
We have $r = \dim \ker T(\lambda_0) = 2$. For a general analytic vector function $y(\lambda)$ with power expansion 
$y^0 + y^1 \lambda + y^2 \lambda^2 + \cdots$ around $0$, we can write (specifying explicitly only the terms up to order 
$\lambda^2$)
$$
T(\lambda)y(\lambda) = 
\begin{bmatrix}
(y^0_1 + y^0_2) \lambda^2  \\ 
(y^0_2 + y^0_3) \lambda + (y^0_1 + y^1_2 + y^1_3) \lambda^2 \\
y^0_3 + (y^0_1 + y^1_3) \lambda + (y^1_1 + y^2_3) \lambda^2 
\end{bmatrix} + \cdots\,.
$$
The successive equations to be fulfilled for $Ty$ to have a zero of order 1, 2, 3 at $\lambda=0$ are:\footnote{%
\nw{That is to say, (\ref{eqn1}) is the equation to be fulfilled for a zero of order 1, (\ref{eqn1}) and (\ref{eqn2}) together are the
equations for a zero of order 2, and so on.}}
\begin{subequations}
\begin{align}
& y^0_3 = 0 \label{eqn1} \\
& y^0_2 + y^0_3 = 0, \quad  y^0_1 + y^1_3 = 0 \label{eqn2} \\
& y^0_1 + y^0_2 = 0, \quad y^0_1 + y^1_2 + y^1_3 = 0, \quad y^1_1 + y^2_3 = 0. \label{eqn3}
\end{align}
\end{subequations}
Determining the dimensions of the subspaces $L_j$ comes down to counting the number of independent restrictions that these 
equations impose on the vector $y^0$. Equation (\ref{eqn1}) gives rise to one restriction, (\ref{eqn2}) generates
another one, and (\ref{eqn3}) adds a third. Therefore, we find $\ell_1 = 2$, $\ell_2 = 1$, and $\ell_3 = 0$. From the rule (\ref{ki}), we 
therefore have $m_1 = 2$, $m_2 = 1$, and $m_3 = 0$.
\end{example}

For any matrix $M\in \cH_0^{n \times m}$, one
can unambiguously define an induced mapping $\overline M$ from the quotient space $\cM^m/\cH_0^m$ to the
quotient space $\cM^n/\cH_0^n$ by\footnote{\nw{The equivalence class of $y \in \cM^n$ modulo $\cH_0^n$ is
denoted by $[y]$. If one prefers to think in term of representatives rather than equivalence classes, the notation 
$[y]$ can also be read as ``principal part of $y$ at $\lambda_0$''\!.}}
\begin{equation} \label{T0def}
\overline M \cn [y] \mapsto [My] \qquad\qquad (y \in \cM^m).
\end{equation} 
This is a linear mapping. 
To connect the dimension of $\ker \overline{T}$ to the \nw{partial} multiplicities defined
in (\ref{can}), we need the following simple lemma.

\begin{lemma} \label{inj}
If $M \in \cH_0^{n \times m}$ is left unimodular, then the induced mapping $\overline M \cn \cM^m/\cH_0^m 
\rightarrow \cM^n/\cH_0^n$ is injective.
\end{lemma}
\begin{proof}
To prove the claim, we have to show that, if $y \in \cM^m$ is such that $My \in \cH_0^n$, then $y \in \cH_0^m$. But 
this is immediate, since $M$ has a left inverse in $\cH_0^{m \times n}$.
\end{proof}


\begin{proposition}  \label{dimkerT}
We have
\begin{equation} \label{dimker}
\dim \ker \oT = \sum_{i=1}^r m_i.
\end{equation} 
\end{proposition}
\begin{proof}
If $U \in \cH_0^{n \times n}$ is unimodular, then, for any $[y] \in \ker \oT$, we have $[U^{-1}y] \in \ker \overline{TU}$
and $\overline{U}[U^{-1}y] = [y]$. This shows that the induced mapping $\overline{U}\cn \cM^n/\cH_0^n \rightarrow
\cM^n/\cH_0^n$ maps the space $\ker \overline{TU}$ onto the space
$\ker \oT$. Since $\overline{U}$ is injective by the above lemma, it follows that $\dim \ker \overline{TU} =
\dim \ker \oT$. Again by the lemma above, we also have $\ker \overline{UT} = \ker \oT$ when $U$ is 
unimodular. We can therefore use the Smith form to conclude that (\ref{dimker}) holds.
\end{proof}

The finite-dimensional space $\ker \oT$ is used in the following representation of the principal part of $T^{-1}$ at $\lambda_0$ 
in terms of linear mappings.
\begin{proposition} \label{ABC}
Write $X = \ker \oT$.
Let mappings $A \cn X \rightarrow X$, $B \cn \mC^n \rightarrow X$, and 
$C \cn X \rightarrow \mC^n$ be defined by  
\begin{equation} \label{ABCdef}
A \cn [y] \mapsto [\chi_0 y], \qquad B\cn u \mapsto [T^{-1}u], \qquad C \cn [y] = y_{-1}
\end{equation}
where $y_{-1}$ stands for the coefficient of $\chi_0^{-1}$ in the Laurent series expansion of $y \in \cM^n$
around $\lambda_0$. The following representation holds:
\begin{equation} \label{statespace}
T^{-1} \doteq C(\chi_0 I - A)^{-1}B.
\end{equation}
\end{proposition}
\begin{proof}
Note that $A$ is well-defined since multiplication by $\chi_0$ maps $\cH_0$ into itself, $B$ indeed maps into $\ker \oT$ since
$T(T^{-1}u) = u \in \cH_0^n$ for all $u \in \mC^n$, and $C$ is well-defined since knowledge of only the principal part of $y$
at $\lambda_0$ is enough to determine $y_{-1}$. Let $s$ denote the pole order of $T^{-1}$ at $\lambda_0$, and 
write the principal part of $T^{-1}$ at $\lambda_0$ as $\sum_{k=1}^s R_k \chi_0^{-k}$. For $u \in \mC^n$, we have
$$
CA^{j-1}Bu = (\chi_0^{j-1}T^{-1}u)_{-1} = R_j u \qquad (j=1,2,\dots).
$$
If $y \in \cM^n$ is such that $Ty \in \cH_0$, then $\chi_0^s y = \chi_0^s T^{-1} Ty \in \cH_0^n$; hence, $A^s=0$.
Consequently, 
$$
C(\chi_0 I - A)^{-1}B = C \Big( \sum_{j=1}^s A^{j-1}\chi_0^{-j} \Big) B = \sum_{j=1}^s  CA^{j-1}B \chi_0^{-j} 
= \sum_{j=1}^s  R_j \chi_0^{-j} \doteq T^{-1}.
$$
\vskip-4mm
\end{proof}
\noindent
In the language of mathematical system theory, the representation (\ref{statespace}) is a \emph{state space realization}
of the principal part of $T^{-1}$. The proposition above uses the realization method of \cite{Fuhrmann1976,Fuhrmann2002},
together with the observation of \cite{KS} that Fuhrmann's method works just as well when the rings of polynomials and
strictly proper rational functions that he used are replaced by the ring $\cH_0$ and the ring of
polynomials in $\chi_0^{-1}$ with zero constant term.\footnote{\nw{\cite{KS} actually work with rational (rather than meromorphic)
functions, but this is inessential.}} 

To work out the result of Prop.\,\ref{ABC} in a concrete case such as (\ref{Tdef}), one needs to choose a suitable basis in the
vector space $\ker \oT$, preferably in such a way that matrix representations of the mappings $A$, $B$, and $C$ can be computed
using information from the holomorphic matrix $T$ directly, rather than by computing its inverse as might be suggested
in particular by the definition of $B$ in (\ref{ABCdef}). This is in fact one of the most important applications of Keldysh's 
theorem, as discussed below. An example to illustrate the proposition is therefore deferred until later (Example \ref{xmpl4}).

\section{Root functions and canonical matrices}

The following definition is standard in the literature; see for instance \cite[Def.\,1.6]{MM}.%
\footnote{One might use $\cH_0$ rather than $\cH$ in the definition of root 
functions without causing harm. See also Lemma \ref{Hn}.}

\begin{definition} \label{rootfunction}
A function $y \in \cH^n$ is said to be a \emph{root function} for $T$ at $\lambda_0$ if $y(\lambda_0) \neq 0$
and $T(\lambda_0)y(\lambda_0) = 0$. The order of the zero of $Ty$ at $\lambda_0$
is called the \emph{multiplicity} of $y$ (with respect to $T$ at $\lambda_0$), and is denoted by $\nu(y)$.
\end{definition}
\noindent
Gohberg and Sigal \cite{GS} attribute the introduction of root functions to Kre{\u\i}n and Trofimov \cite{Krein} and 
Macaev and Palant \cite{Palant} independently.\footnote{%
The latter reference is as given in \cite{GS}, but is authored only by Palant. Perhaps, Gohberg and Sigal had \cite{Macaev} in mind.
The source cited by Mennicken and M\"{o}ller \cite{MM} for the notion of root functions is Trofimov \cite{Trofimov}.} 
Keldysh \cite{Keldysh1951,Keldysh1971} 
instead uses ``systems of eigenvectors and associated vectors''\!. Such a system appears as the sequence of 
coefficients of a polynomial root function.\footnote{Requiring root functions to be polynomial
is possible (see the proof of Lemma \ref{Hn}); this restriction is often viewed as inconvenient, however.} 
The coefficient corresponding to the constant term of the polynomial
is called an eigenvector, since it must belong to $\ker T(\lambda_0)$, and the other coefficients are 
called associated vectors. In part of the literature (for instance \cite{BGK}), a system of eigenvectors and 
associated vectors is called a \emph{Jordan chain}. 

In Def.\,\ref{rootfunction}, root functions are required to be holomorphic throughout the domain $\Omega$,
whereas, in the definition of the subspaces $L_j$ in (\ref{LjT}), locally holomorphic functions were used. The
following lemma (which formalizes a remark in \cite[p.\,14]{MM}) shows that this difference is inconsequential.

\begin{lemma} \label{Hn}
A vector $y^0 \in \mC^n$ with $y^0 \neq 0$ belongs to the subspace $L_j$ if and only if there exists 
a root function $y \in \cH^n$ of multiplicity $j$ such that $y(\lambda_0) = y^0$.
\end{lemma}

\begin{proof}
If $y^0 \in L_j$, then, by definition, there exists $\ty \in \cH_0^n$ such that $\ty(\lambda_0)=y^0$ and 
$T\ty \in \chi_0^j\cH_0^n$. For a function 
$y \in \cH^n$ to be such that $y(\lambda_0)=y_0$ and $Ty \in \chi_0^j\cH_0^n$, it is sufficient that 
$y^{(k)}(\lambda_0) = \ty^{(k)}(\lambda_0)$ for $k = 0,1,\dots,j$. These requirements can indeed be satisfied; 
in particular, one can choose $y$ to be a polynomial. This proves the ``only if'' part of the claim; the ``if'' part is trivial.
\end{proof}

\nw{To make Prop.\,\ref{ABC} operational, collections of root functions with special properties are required. A first observation
is as follows.}
\begin{proposition} \label{ub}
If $\{y_1,\dots,y_q\}$ is a collection of root functions for $T$ at $\lambda_0$ such that the constant vectors
$\{y_1(\lambda_0),\dots,y_q(\lambda_0)\}$ are linearly independent, then the sum of the multiplicities of
these root functions does not exceed the algebraic multiplicity of the root of $T$ at $\lambda_0$.
\end{proposition}

\begin{proof}
Let $Y \in \cH_0^{n\times q}$ be defined as the matrix with columns $y_i$ ($i=1,\dots,q$), write $\nu_i = \nu(y_i)$,
and define $D = \diag(\chi_0^{\nu_1},\dots,\chi_0^{\nu_q})$. Note that $Y$ is left unimodular. By the definition of
root functions, we have $TY\!D^{-1} \in \cH_0^{n \times q}$, which implies that the space 
$(Y\!D^{-1}\cH_0^q)/\cH_0^n \subset
\cM^n/\cH_0^n$ is a subspace of $\ker \overline{T}$. Since the induced mapping $\overline{Y}$ is injective on
$\cM^q/\cH_0^q$ by Lemma \ref{inj}, the dimension of $(Y\!D^{-1}\cH_0^q)/\cH_0^n$ is the same as the dimension
of $(D^{-1}\cH_0^q)/\cH_0^q$. The matrix $D$ induces a bijective mapping from 
$(D^{-1}\cH_0^q)/\cH_0^q$ to $\cH_0^n/(D\cH_0^n)$. The dimension of the latter space is equal to $\sum_{i=1}^q
\nu_i$, since it has a basis consisting of equivalence classes modulo $D\cH_0^n$ of the form $[\chi_0^j e_i]$, where 
$e_i$ denotes the $i$-th unit vector in $\cH_0^q$ ($i=1,\dots,q$) and $0 \leq j \leq \nu_i-1$. It follows that 
$$
\sum_{i=1}^q \nu_i = \dim\, (D^{-1}\cH_0^q)/\cH_0^q = \dim\, (Y\!D^{-1}\cH_0^q)/\cH_0^q \leq 
\dim \ker \overline{T} = \sum_{i=1}^r m_i.
$$
\vskip-4mm
\end{proof}

\noindent
The proposition \nw{is part of the motivation for} the following definition.
\begin{definition}
A \emph{canonical system of root functions} for $T$
at $\lambda_0$ is a collection $\{y_1,\dots,y_r\}$ of root functions for $T$ at $\lambda_0$ such that 
the following conditions are satisfied:
\begin{itemize}
\item[(i)] the vectors $y_1(\lambda_0),\dots, y_r(\lambda_0) \in \mC^n$ are linearly independent; 
\item[(ii)] $\sum_{i=1}^r \nu(y_i) = \sum_{i=1}^r m_i$;
\item[(iii)] the multiplicities $\nu(y_i)$ are in nonincreasing order, i.e., $\nu(y_i) \geq \nu(y_k)$ whenever $i \leq k$.
\end{itemize}
\end{definition}
\vskip2mm\noindent
Condition (iii) can of course always be satisfied by reordering the root functions if necessary; the
requirement just serves to simplify the notation.

We shall say that a basis $\{y_1^0, \dots, y_r^0\}$ of $\ker T(\lambda_0)$ is \emph{adapted} to the nonincreasing 
sequence $\{L_j\}_{1 \leq j \leq s}$ when, for each $j=1,\dots,s$, the subspace $L_j$ is spanned by the vectors
$y_1^0,\dots,y_{\ell_j}^0$ where $\ell_j = \dim L_j$. Clearly, one can construct such a basis by starting with a basis 
for the smallest subspace $L_s$, then extending it to a basis for the next larger subspace $L_{s-1}$, and so on.

\begin{proposition} \label{cansys}
A canonical system of root functions for $T$ at $\lambda_0$ exists. 
\end{proposition}

\begin{proof}
Take a basis $\{ y_1^0, \dots, y_r^0 \}$ for $\ker T(\lambda_0)$ that is adapted to the
sequence $L_1,\dots,L_s $. By Lemma \ref{Hn}, there are root functions $y_1,\dots,y_r$ such that $y_i(\lambda_0)
= y_i^0$ for $i=1,\dots,r$ and $\nu(y_i) \geq j$ when $y_i^0 \in L_j$. Because of the adaptedness of the basis,
the root functions that are constructed in this way contain $\ell_{m_1}$ elements with multiplicity at least $m_1 = s$,
$\ell_{m_2}-\ell_{m_1}$ elements with multiplicity at least $m_2$, and so on. Therefore, we have (using (\ref{summi}))
\begin{equation} \label{sumineq}
\sum_{i=1}^r \nu(y_i) \geq \sum_{j=1}^s j(\ell_j-\ell_{j+1}) = \sum_{j=1}^s \ell_j = \sum_{i=1}^r m_i.
\end{equation}
By Prop.\,\ref{ub}, equality holds, and the system $\{y_1,\dots,y_r\}$ is canonical.
\end{proof}

\begin{proposition} \label{canchar}
If $\{y_1,\dots,y_r\}$
is a canonical system, then  
$\nu(y_i) = m_i$ for $i=1,\dots,r$.
\end{proposition}
\begin{proof}
Write $\nu_i := \nu(y_i)$. 
For $j=1,\dots,s$, define $\hl_j = \#\{ 1 \leq i \leq r \mid \nu_i \geq j \}$. The sequences $(\nu_1,\dots,\nu_r)$ and
$(\hl_1,\dots\hl_s)$ are conjugate partitions, just as $(m_1,\dots,m_r)$ and $(\ell_1,\dots,\ell_s)$ are. Moreover,
since $y_i(\lambda_0) \in L_j$ when $\nu_i \geq j$, we have $\hl_j \leq \ell_j$ for all $j=1,\dots,s$.
It follows that
$$
\nu_i = \#\{ 1 \leq j \leq s \mid \hl_j \geq i \} \,\leq\, \#\{ 1 \leq j \leq s \mid \ell_j \geq i \} = m_i
\qquad\qquad (i=1,\dots,r).
$$    
Consequently, the canonicity condition $\sum_{i=1}^r \nu_i = \sum_{i=1}^r m_i$ can hold only when
$\nu_i = m_i$ for all $i=1,\dots,r$.
\end{proof}

It is convenient to employ a matrix format for canonical systems of root functions. 

\begin{definition} \label{canmatdef}
A matrix $Y \in \cH^{n \times r}$ is said to be a \emph{right canonical matrix} (for $T$ at $\lambda_0$) when its
columns form a canonical system of root functions for $T$ at $\lambda_0$.
\end{definition}
\noindent
The existence of right canonical matrices follows from Prop.\,\ref{cansys}.

\begin{example} \label{xmpl2}
In the setting of Example \ref{xmpl}, a matrix $Y \in \cH^{3 \times 2}$ is a right canonical matrix when its
first column satisfies equations (\ref{eqn1}) and (\ref{eqn2}), its second column satisfies (\ref{eqn1}), and 
$Y(\lambda_0)$ has full column rank. Therefore, a right canonical matrix is given, for instance, by
\begin{equation} \label{Ydef}
Y(\lambda) = \begin{bmatrix} 1 & 0 \\ 0 & 1 \\ -\lambda & 0 \end{bmatrix}.
\end{equation}
We have
$$
T(\lambda)Y(\lambda) = \begin{bmatrix} \lambda^2 & \lambda^2 & 0 \\ \lambda^2 & \lambda & \lambda \\ \lambda & 0 & 1
\end{bmatrix} \begin{bmatrix} 1 & 0 \\ 0 & 1 \\ -\lambda & 0 \end{bmatrix} = 
\begin{bmatrix} \lambda^2 & \lambda^2 \\ 0 & \lambda \\ 0 & 0 \end{bmatrix}.
$$
The first column is divisible (in $\cH_0$) by $\lambda^2$ and the second by $\lambda$, as required.
\end{example}

Since the root functions in a canonical system are ordered by nonincreasing multiplicity, we have the following.

\begin{lemma} \label{final}
If $Y$ is a right canonical matrix for $T$ at $\lambda_0$, then, for each $j=1,\dots,s$, the first $\ell_j$ columns 
of $Y(\lambda_0)$ form a basis for the subspace $L_j$.
\end{lemma}
\noindent
Recall the definition of the matrix $\Delta$ in (\ref{deltadef}) and the notation $\doteq$ introduced in (\ref{doteq}).
\begin{proposition} \label{canmat}
A matrix $Y \in \cH^{n \times r}$ is a right canonical matrix for $T$ at $\lambda_0$ if and only if the
following two conditions are satisfied:
\begin{itemize}
\item[(i)]
$Y(\lambda_0) \in \mC^{n \times r}$ has full column rank;
\item[(ii)] $TY\!\Delta^{-1} \doteq 0$.
\end{itemize}
\end{proposition}

\begin{proof}
Suppose that conditions (i) and (ii) hold.
Condition (ii) is equivalent to the statement that there exists a matrix $H \in \cH_0^{n \times r}$ such
that $TY=H\Delta$; i.e., $Ty_i=h_i\chi_0^{m_i}$ for $i = 1,\dots, r$, where $y_i \in \cH^n$ and $h_i \in \cH_0^n$ 
denote the $i$-th columns of $Y$ and $H$ respectively. This implies that $y_i$ is a root function for $T$ at $\lambda_0$, with multiplicity $\nu(y_i) \geq m_i$. By condition (i) and Prop.\,\ref{ub}, the sum of the multiplicities of the 
root functions $y_i$ is equal to $\sum_{i=1}^r m_i$, and hence $\{y_1,\dots,y_r\}$ is a canonical system. 
Conversely, if $Y \in \cH^{n\times r}$ is right canonical, then, for each $i=1,\dots,r$, the $i$-th column of $Y$ satisfies
$Ty_i \in \chi_0^{m_i}\cH_0^n$ by Prop.\,\ref{canchar}. This implies that (ii) holds, and (i) holds by definition.
\end{proof}
\noindent
The following proposition rephrases and extends Lemma A.9.3 of \cite{Kozlov} in transposed form. The proof
uses a modification of the argument in the proof of Thm.\,1.2 in \cite{GKvS}.

\begin{proposition} \label{crux}
Let $Y$ be a right canonical matrix for $T$ at $\lambda_0$. Then the constant matrix
$(TY\!\Delta^{-1})(\lambda_0) \in \mC^{n\times r}$ has full column rank, and  the following direct sum 
decomposition holds:
\begin{equation}
\mC^n = \im T(\lambda_0) \oplus \im\kern1pt (TY\!\Delta^{-1})(\lambda_0).
\end{equation} 
\end{proposition}
\begin{proof}
For a proof by contradiction,
suppose there exist $y^0 \in \mC^n$ and $x \in \mC^r$, with $x \neq 0$, such that 
$(TY\!\Delta^{-1})(\lambda_0) x = T(\lambda_0) y^0$. Define $k = \min \{ j \geq 1 \mid
\chi_0^j \Delta^{-1}x \in \cH_0^r \}$. From $\chi_0^k \Delta^{-1} x \in \cH_0^n$ it follows that $x_i = 0$ for 
all indices $i$ such that $m_i \geq k+1$, i.e., $1 \leq i \leq \ell_{k+1}$.  
Write $\hx := (\chi_0^k \Delta^{-1} x)(\lambda_0)$. Note that $\hx \neq 0$, because otherwise we would
have $\chi_0^{k-1}\Delta^{-1}x \in \cH_0^r$ and $k$ would not be minimal. 
Define $y = \chi_0^k(Y\!\Delta^{-1}x-y^0) \in \cH_0^n$.
The assumption $(TY \Delta^{-1})(\lambda_0) x = T(\lambda_0) y^0$ implies that $TY\!\Delta^{-1}x - Ty^0
\in \chi_0 \cH_0^n$, and hence $Ty = \chi_0^k(TY\!\Delta^{-1}x - Ty^0) \in \chi_0^{k+1}\cH_0^n$.
It follows that $y(\lambda_0) \in L_{k+1}$. On the other hand, $y(\lambda_0) = Y(\lambda_0)\hx$, and 
the vector $\hx$ inherits from $x$ the property that all of its entries with indices $i \leq \ell_{k+1}$ are 0.
In view of Lemma \ref{final}, we have a contradiction. It follows that the matrix $(TY\!\Delta^{-1})(\lambda_0)$
has full column rank (for this, take $y^0=0$) and that 
$\im T(\lambda_0) \cap \im\kern1pt (TY\!\Delta^{-1})(\lambda_0) = \{ 0 \}$. 
The proof is completed by noting that the dimension of $\im T(\lambda_0)$ is $n-r$,
while $\dim \im\kern1pt (TY\!\Delta^{-1})(\lambda_0) = r$ by the full column rank property.
\end{proof}

The following proposition shows how canonical matrices give rise to a basis for the space $\ker \oT$. It is
shown further below (see Remark \ref{ABCrep} and Example \ref{xmpl4}) that the
mappings $A$, $B$, and $C$ appearing in the abstract state space representation (\ref{statespace}) can
conveniently be expressed in matrix form with respect to this basis. 
\begin{proposition} \label{basis}
Let $Y$ be a right canonical matrix for $T$ at $\lambda_0$. Then a basis over $\mC$ for the space $\ker \oT$
is given by the elements $[Y\Delta^{-1}\chi_0^{j-1} e_i] \in \cM^n/\cH_0^n$ 
($i=1,\dots,r$, $j=1,\dots,m_i$), where $e_i$ denotes the $i$-th unit vector in $\mC^r$.
\end{proposition}

\begin{proof}
The number of elements $[Y\Delta^{-1}\chi_0^{j-1} e_i]$ ($i=1,\dots,r$, $j=1,\dots,m_i$) is $\sum_{i=1}^r m_i$,
which is equal to $\dim\ker \oT$ by Prop.\,\ref{dimkerT}. It remains to show that these elements are linearly independent.
To this end, take arbitrary $v_i^j \in \mC$ ($i=1,\dots,r$, $j=1,\dots,m_i$), and define $v \in \cH_0$ by $v = 
\sum_{i=1}^r \sum_{j=1}^{m_i} \chi_0^{j-1} e_i v_i^j$. Suppose that $[Y\Delta^{-1}v]=0$, i.e., $Y\Delta^{-1}v \in \cH_0^n$.
Note that $\Delta^{-1}v$ is a polynomial vector in $\chi_0^{-1}$
with zero constant term. On the other hand, since $Y$ is left unimodular by Prop.\,\ref{canmat}(i), the
assumption $Y\Delta^{-1}v \in \cH_0^n$ implies that $\Delta^{-1}v \in \cH_0^r$. Consequently, $\Delta^{-1}v=0$,
and hence $v=0$. Since the elements $\chi_0^{j-1}e_i$ are linearly 
independent, it follows that $v_i^j = 0$ for all $i=1,\dots,r$ and $j=1,\dots,m_i$. This proves the desired result.
\end{proof}

\section{Keldysh's theorem}

A holomorphic row vector $v$ of length $n$ is said to be a \emph{left
root function} for $T$ at $\lambda$ when $v(\lambda_0) \neq 0$ and $v(\lambda_0)T(\lambda_0)=0$. The notion of a 
left canonical matrix for $T$ at $\lambda_0$ is defined analogously to Def.\,\ref{canmatdef}, and one has the
following characterization.

\begin{proposition} \label{leftcan}
A matrix $V \in \cH^{r \times n}$ is a left canonical matrix for $T$ at $\lambda_0$ if and only if 
$V(\lambda_0) \in \mC^{r \times n}$ has full row rank, and
$\Delta^{-1}VT \doteq 0$.
\end{proposition}
\noindent
Prop.\,\ref{crux} translates to left canonical matrices as follows.

\begin{proposition} \label{crux2}
Let $V$ be a left canonical matrix for $T$ at $\lambda_0$. Then the constant matrix $(\Delta^{-1}VT)(\lambda_0)$ has full row rank, and the following direct sum decomposition holds:
\begin{equation} \label{dirsum}
\mC^n = \ker T(\lambda_0) \oplus \ker\kern1pt (\Delta^{-1}VT)(\lambda_0).
\end{equation}
\end{proposition}
\vskip2mm\noindent
For later use, we note the following corollary.
\begin{corollary} \label{VTYD}
Let $V$ and $Y$ be left and right canonical matrices, respectively, for $T$ at $\lambda_0$. Then the
matrices $\Delta^{-1}VTY \in \cH_0^{r \times r}$ and $VTY\!\Delta^{-1} \in \cH_0^{r \times r}$ are unimodular.
\end{corollary}
\begin{proof}
From the definition of right canonical matrices, it follows that $Y(\lambda_0) \in \mC^{n \times r}$ is a basis
matrix for $\ker T(\lambda_0)$. Since the subspace $\ker\,\Delta^{-1}VT(\lambda_0)$ is complementary to
$\ker T(\lambda_0)$ by Prop.\,\ref{crux2}, the matrix $(\Delta^{-1}VTY)(\lambda_0) = (\Delta^{-1}VT)(\lambda_0)
Y(\lambda_0) \in \mC^{r \times r}$ is nonsingular. Consequently, $\Delta^{-1}VTY$ is unimodular. The proof
for $VTY\!\Delta^{-1}$ is analogous.
\end{proof}

\vskip2mm
The Smith form may be used to show that the principal part of $T^{-1}$ at $\lambda_0$ can be described
in terms of suitably chosen left and right canonical matrices for $T$ at $\lambda_0$.

\begin{theorem} \label{gs}
There exist left and right canonical matrices $V$ and $Y$ for $T$ at $\lambda_0$ such that 
$T^{-1} \doteq Y\!\Delta^{-1}V$.
\end{theorem}
\begin{proof}
There are unimodular matrices $U_L, \, U_R \in \cH_0^{n \times n}$ such that, with conformable partitioning,
\begin{equation} \label{gs1}
T = \begin{bmatrix} U_L^1 & U_L^2 \end{bmatrix} \begin{bmatrix} \Delta & 0 \\ 0 & I_{n-r} 
\end{bmatrix} \begin{bmatrix} U_R^1 \\ U_R^2 \end{bmatrix}.
\end{equation}
Write, again with conformable partitioning,
\begin{equation} \label{gs2}
U_L^{-1} = \begin{bmatrix} V \\ \tV \end{bmatrix}, \qquad U_R^{-1} = \begin{bmatrix} Y & \tY \end{bmatrix}.
\end{equation}
We have
\begin{equation} \label{Tinv}
T^{-1} = \begin{bmatrix} Y & \tY \end{bmatrix}  \begin{bmatrix} \Delta^{-1} & 0 \\ 0 & I_{n-r} \end{bmatrix} 
\begin{bmatrix} V \\ \tV \end{bmatrix} = Y\!\Delta^{-1} V + \tY \tV \doteq Y\!\Delta^{-1}V. 
\end{equation}
From the fact that $[Y(\lambda_0) \;\; \tY(\lambda_0)]$ is nonsingular it follows that $Y(\lambda_0)$ has full 
column rank $r$. Multiplying from the left by $T$ and from the right by $U_L^1$,
one finds from (\ref{Tinv}) that $TY\!\Delta^{-1} \doteq 0$,
so that $Y$ is right canonical by Prop.\,\ref{leftcan}. Left canonicity of $V$ is proved analogously.
\end{proof}
\noindent
\begin{remark} \label{VTYw}
For future reference, note that from (\ref{gs1}) and (\ref{gs2}) it also follows that
\begin{equation}
\begin{bmatrix} \Delta & 0 \\ 0 & I_{n-r} \end{bmatrix} =
\begin{bmatrix} V \\ \tV \end{bmatrix} T \begin{bmatrix} Y & \tY \end{bmatrix} =
\begin{bmatrix} VTY & VT\tY \\ \tV TY & \tV T\tY \end{bmatrix}.
\end{equation}
In particular, we have $VTY=\Delta$.
\end{remark}

Thm.\,\ref{gs} is the matrix version of Thm.\,7.1 in \cite{GS}, and the proof as given above follows their argument.
The result given in \cite{Keldysh1951,Keldysh1971} has a stronger matrix version. We write it as follows. 

\begin{theorem}[Keldysh 1951] \label{keldysh}
For any right canonical matrix $Y$ for $T$ at $\lambda_0$, there exists a left canonical matrix $V$ for $T$ at 
$\lambda_0$, determined uniquely up to addition of a matrix of the form $\Delta H$ with $H \in \cH^{r \times n}$, 
such that 
\begin{equation} \label{rpp}
T^{-1} \doteq Y\!\Delta^{-1}V.
\end{equation}
For any such matrix $V$, one has
\begin{equation} \label{biorth}
\Delta^{-1}VTY\!\Delta^{-1} \doteq \Delta^{-1}.
\end{equation}
\end{theorem}
\vskip2mm\noindent
Condition (\ref{biorth}) is called the \emph{biorthogonality condition}. It is stated here in essentially
the form given by Kozlov and Maz'ya \cite[Thm.\,A.10.1]{Kozlov}. These authors also prove (Thm.\,A.10.2) that the condition 
(\ref{biorth}) is
equivalent to the quadruply indexed formulation in terms of systems of eigenvectors and associated vectors as 
given by Keldysh \cite{Keldysh1951,Keldysh1971} and by Mennicken and M\"{o}ller \cite{MM}.

To prove the theorem, one can make use of the following lemma. The proof as given below is a rephrasing of
part of the argument in the proof of Thm.\,1.5.4 in \cite{MM}. 

\begin{lemma} \label{MMlemma}
Let $Y \in \cH^{n \times r}$ be a right canonical matrix for $T$ at $\lambda_0$, and suppose $V \in \cH^{r \times n}$
is such that $T^{-1}-Y\!\Delta^{-1}V$ has pole order $j \geq 1$ at $\lambda_0$. Then there exists a matrix $\tV
\in \cH^{r \times n}$ such that $T^{-1}-Y\!\Delta^{-1}\tV$ has pole order  at most $j-1$ at $\lambda_0$.
\end{lemma}

\begin{proof}
Write $A = T^{-1}-Y\!\Delta^{-1}V$. By assumption, there are matrices $A_1,\dots,A_j \in \mC^{n \times n}$ such that 
$A \doteq A_1 \chi_0^{-1} + \cdots + A_j \chi_0^{-j}$. For any $x \in \mC^n$, we have $A_j x = 
(\chi_0^j Ax)(\lambda_0) \in L_j$, since $T \chi_0^j Ax = \chi_0^j(I-TY\!\Delta^{-1}V)x \in \chi_0^j \cH_0^n$
by Prop.\,\ref{canmat}. This implies, by Lemma \ref{final}, that each of the columns of $A_j$ can be written as a 
linear combination of the first $\ell_j$ columns of $Y(\lambda_0)$. Let the submatrix of $Y$ formed by these columns
be denoted by $Y_{1:\ell_j}$, and let
$M \in \mC^{\ell_j \times n}$ be such that $A_j = Y_{1:\ell_j}(\lambda_0) M$. 
Then the pole order of $A- \chi_0^{-j}Y_{1:\ell_j}M$ at $\lambda_0$ is at most $j-1$. Now define
$$
\tV = V + \chi_0^{-j} \Delta \begin{bmatrix} M \\ 0_{(r-\ell_j)\times n} \end{bmatrix}.
$$
Note that $\tV \in \cH^{r\times n}$, since $m_i \geq j$ for $i \leq \ell_j$. Moreover, 
$$
T^{-1}-Y\!\Delta^{-1}\tV = A - Y\!\Delta^{-1}\chi_0^{-j} \Delta \begin{bmatrix} M \\ 0 \end{bmatrix} =
A - \chi_0^{-j}Y_{1:\ell_j}M 
$$
so that $\tV$ has the required properties.
\end{proof}
\noindent
Keldysh's theorem, in the matrix case, can now be proved as follows.

\begin{proof}
Existence of $V \in \cH^{r\times n}$ such that (\ref{rpp}) holds follows by repeated application of Lemma 
\ref{MMlemma}, starting with $V=0$. For the rest of the proof, note that $Y$ is left unimodular by 
Prop.\,\ref{canmat}; let
$Z \in \cH_0^{r \times n}$ be a left inverse. The biorthogonality condition (\ref{biorth}) is obtained by
multiplying both sides of (\ref{rpp}) from the left by $Z$ and from the right by $TY\!\Delta^{-1} \in \cH_0^{n\times r}$. 
Next, if $V$ satisfies (\ref{rpp}), then clearly the same holds for any matrix $\hV$ of the form $\hV=V+\Delta H$
where $H \in \cH^{r \times n}$. Conversely, suppose that $V,\,\hV \in \cH^{r\times n}$ both satisfy (\ref{rpp}).
Define $H = \Delta^{-1}(\hV-V)$; then $\hV = V+\Delta H$. We have $Y\!H \doteq 0$; multiplication from the left by 
$Z$ shows that $H$ belongs to $\cH_0^{r \times n}$. Since $H$ is of the form $H = \Delta^{-1}\tV$ with
$\tV \in \cH^{r \times n}$, and $\Delta$ only has zeros at $\lambda_0$, it follows that in fact $H \in \cH^{r \times n}$.

By multiplying both sides of (\ref{rpp}) from the left by $Z$ and from the right by
$T$, one finds that $\Delta^{-1}VT \doteq 0$. Furthermore, the biorthogonality relation (\ref{biorth}) means that
there exists a matrix $H \in \cH_0 ^{r \times r}$ such that $\Delta^{-1} = \Delta^{-1}VTY\!\Delta^{-1} + H$, which
implies that $VTY\!\Delta^{-1} = I_r-\Delta H$. Since $\Delta(\lambda_0)=0$ and 
$TY\!\Delta^{-1} \in \cH_0^{n \times r}$, it follows that $V(\lambda_0) (TY\!\Delta^{-1})(\lambda_0) = I_r$,
which implies that $V(\lambda_0)$ has full row rank. Canonicity of $V$ now follows from Prop.\,\ref{leftcan}. 
\end{proof}

\begin{remark} \label{ABCrep}
Using the basis constructed in Prop.\,\ref{basis}, matrix versions of the abstract first-order representation of
Prop.\,\ref{ABC} can be obtained from (\ref{rpp}). In particular, the representation (\ref{rpp}) allows the mapping $B$
appearing in (\ref{ABCdef}) to be written as $B\cn  u \mapsto [Y\Delta^{-1}Vu]$, so that the entries of a matrix
representation of $B$ can be found by expressing the columns of $V$ as linear combinations of the linearly independent
vectors $\chi_0^{j-1}e_i$ ($i=1,\dots,r$, $j=1,\dots,m_i$). Note that one can indeed do this after possibly replacing
$V$ by a matrix $V+\Delta H$ with $H \in \cH_0^{r \times n}$. Such a modification is admissible since the canonicity 
property is preserved (this follows from Prop.\,\ref{canmat}), and moreover the representation (\ref{rpp}) remains valid. 
An example of the construction of a matrix representation is given below (Example \ref{xmpl4}).
\end{remark}

\section{Extensions}

Kozlov and Maz'ya \cite[Thm.\,A.10.1]{Kozlov} present Keldysh's theorem with an additional statement concerning the uniqueness
of solutions of the biorthogonality condition (\ref{biorth}), read as an equation for the left canonical matrix $V$ when
the right canonical matrix $Y$ is given. They restrict the choice of left and right canonical matrices
$V$ and $Y$ to the ones for which $\Delta^{-1}V$ and $Y\!\Delta^{-1}$ are polynomial matrices in $\chi_0^{-1}$
\nw{with zero constant term}.
This restriction is avoided in the version below. Also, the proof provides more detail than is given in \cite{Kozlov}.
For notational convenience, we look at the biorthogonality condition as an equation for the right canonical matrix
when the left canonical matrix is given. The following simple observation will be used.

\begin{lemma} \label{obs}
Let $M \in \cH_0^{p \times n}$ be left unimodular, and let $j \geq 1$. If $y \in \cH_0^n$ 
satisfies $My \in \chi_0^j \cH_0^p$, then $y \in \chi_0^j \cH_0^n$.
\end{lemma}

\begin{proof}
Let $\tM \in \cH_0^{n \times p}$ be a left inverse of $M$. For
$y \in \cH_0^n$ such that $My \in \chi_0^j \cH_0^p$, we have $y = \tM (My) \in \chi_0^j \cH_0^n$.
\end{proof}

\begin{theorem} \label{km}
For any left canonical matrix V for $T$ at $\lambda_0$, there exists a right canonical matrix $Y$ for $T$ at $\lambda_0$
such that the biorthogonality condition (\ref{biorth}) holds. This matrix is determined uniquely up to addition of a 
matrix of the form $H\Delta$ with $H \in \cH^{n \times r}$. Formula 
(\ref{rpp}) for the principal part of $T^{-1}$ at $\lambda_0$ holds with any such matrix.
\end{theorem}

\begin{proof}
Existence has already been shown in Thm.\,\ref{keldysh}. Relation 
(\ref{biorth}) is equivalent to existence of a matrix $\tH \in \cH_0^{r \times r}$ such that
$\Delta^{-1}VTY = I + \tH\Delta$.
To prove uniqueness of right canonical matrices $Y$ that satisfy this equation, up to an additive term of the 
form $H\Delta$ with $H \in \cH^{n \times r}$,
we need to show that, if (i) $\Delta^{-1}VTY = \tH\Delta$ where $\tH \in \cH_0^{r \times r}$ and (ii) $Y$ is the
difference of two right canonical matrices, then $Y = H\Delta$ for some $H \in \cH^{n \times r}$. Condition (i)
can be restated in columnwise fashion as
$$
\Delta^{-1}VT y_i \in \chi_0^{m_i} \cH_0^r \qquad\qquad (i=1,\dots,r)
$$
where $y_i$ denotes the $i$-th column of $Y$. Condition (ii) implies that
$Ty_i \in \chi_0^{m_i} \cH_0^n$ for all $i=1,\dots,r$. By Prop.\,\ref{crux2} and Lemma \ref{obs},
combination of the two conditions leads to the conclusion that $y_i \in \chi_0^{m_i} \cH_0^n$ for all
$i = 1,\dots,r$. In other words,
$Y\!\Delta^{-1}$ does not have a singularity at $\lambda_0$. Since $\Delta$ only has zeros at $\lambda_0$ and
$Y \in \cH^{n \times r}$, it follows that in fact $Y\!\Delta^{-1}$ has no singularities at all in the domain $\Omega$;
i.e., $Y=H\Delta$ for some $H \in \cH^{n \times r}$. 

Finally, let $Y$ be a right canonical matrix for which the biorthogonality condition (\ref{biorth}) 
holds. From the version of Thm.\,\ref{keldysh} in which the left (rather than the right) canonical matrix is taken
to be given, it follows that there exists a right canonical matrix $\tY$ that satisfies the relation 
$T^{-1} \doteq \tY\!\Delta^{-1}V$ as well as the biorthogonality condition. By what has been proved in the previous 
paragraph, there exists $H \in \cH^{n\times n}$ such
that $Y = \tY + H\Delta$. From the fact that $T^{-1} \doteq \tY\!\Delta^{-1}V$, it then follows that we also have
$T^{-1} \doteq Y\!\Delta^{-1}V$.
\end{proof}

\begin{example} \label{xmpl3}
In the setting of Example \ref{xmpl}, one can write equations for left canonical matrices similarly to what was done
in Example \ref{xmpl2} for right canonical matrices. From these, one finds that a left canonical matrix is given,
for instance, by
\begin{equation} \label{Vdef}
V = \begin{bmatrix} 1 & 0 & 0 \\ 0 & 1 & 0 \end{bmatrix}. 
\end{equation}
Using this choice, one finds for a general matrix $Y \in \cH_0^{3 \times 2}$:
$$
(\Delta^{-1}VTY\Delta^{-1})(\lambda) =
\begin{bmatrix}
\lambda^{-2} (y_{11} + y_{21})(\lambda) & \lambda^{-1} (y_{12} + y_{22})(\lambda) \\
\lambda^{-2} (y_{21} + y_{31})(\lambda) + \lambda^{-1} y_{11}(\lambda) & \lambda^{-1} 
(y_{22} + y_{32})(\lambda) + y_{12}(\lambda)
\end{bmatrix}.
$$
The conditions for $Y$ to be right canonical and to satisfy the biorthogonality condition with respect to $V$ are
\begin{subequations}
\begin{align}
& y^0_{11} + y^0_{21} = 1, \quad y^0_{21} + y^0_{31} = 0 \label{e1} \\
& y^1_{11} + y^1_{21} = 0, \quad y^0_{12} + y^0_{22} = 0, \quad y^0_{11} + y^1_{21} + y^1_{31} = 0, \quad 
y^0_{22} + y^0_{32} = 1 \label{e2} \\
& y^0_{31} = 0, \quad y^0_{11} + y^1_{31} = 0, \quad y^0_{21} + y^0_{31} = 0 \label{e3} \\
& y^0_{32} = 0 \label{e4}
\end{align}
\end{subequations}
together with the requirement that $Y(\lambda_0)$ should be of full column rank. Here, equations (\ref{e1}) and 
(\ref{e2}) are due to the biorthogonality condition, (\ref{e3}) corresponds to conditions (\ref{eqn1}) and (\ref{eqn2})
for the first column of $Y$, and (\ref{e4}) is condition (\ref{eqn1}) for the second column of $Y$. The equations lead
to a unique solution for the unknowns that appear in them, namely
$$
Y^0 = \begin{bmatrix} 1 & -1 \\ 0 & 1 \\ 0 & 0 \end{bmatrix}, \qquad
Y^1_{\cdot 1} = \begin{bmatrix} 0 \\ 0 \\ -1 \end{bmatrix}.
$$
As a particular solution, one can therefore take
\begin{equation} \label{bYdef}
\bY(\lambda) = \begin{bmatrix}  1 & -1 \\ 0 & 1 \\ -\lambda & 0 \end{bmatrix}
\end{equation} 
where a bar has been added to distinguish this matrix from the one in (\ref{Ydef}). From Thm.\,\ref{km}, the
principal part of $T^{-1}$ at $\lambda = 0$ is now determined by
\begin{equation} \label{pp}
T^{-1}(\lambda) \doteq \begin{bmatrix}  1 & -1 \\ 0 & 1 \\ -\lambda & 0 \end{bmatrix}
\begin{bmatrix} \lambda^{-2} & 0 \\ 0 & \lambda^{-1} \end{bmatrix} 
\begin{bmatrix} 1 & 0 & 0 \\ 0 & 1 & 0 \end{bmatrix} =
\begin{bmatrix} \lambda^{-2} & - \lambda^{-1} & 0 \\ 0 & \lambda^{-1} & 0 \\ -\lambda^{-1} & 0 & 0 
\end{bmatrix}.
\end{equation}
\end{example}
\vskip2mm
\begin{example} \label{xmpl4}
To illustrate the construction of a representation as in (\ref{statespace}), consider again the holomorphic matrix $T$ of 
Example \ref{xmpl}, and take $V$ in (\ref{Vdef}) and $\bY$ in
(\ref{bYdef}) as left and right canonical matrices satisfying the biorthogonality condition. As a basis for 
$\ker \oT$, we can take the columns (modulo $\cH_0^n$) of the matrix given by
$$
\bY\Delta^{-1} \begin{bmatrix} 1 & \lambda & 0 \\ 0 & 0 & 1 \end{bmatrix} = 
\begin{bmatrix} 1 & -1 \\ 0 & 1 \\ -\lambda & 0 \end{bmatrix} 
\begin{bmatrix} \lambda^{-2} & \lambda^{-1} & 0 \\ 0 & 0 & \lambda^{-1} \end{bmatrix} =
\begin{bmatrix} \lambda^{-2} & \lambda^{-1} & -\lambda^{-1} \\ 0 & 0 & \lambda^{-1} \\
-\lambda^{-1} & -1 & 0 \end{bmatrix}.
$$
Using Remark \ref{ABCrep}, the matrix representations of the mappings $A$, $B$, and $C$ appearing in
Prop.\,\ref{ABC} with respect to this choice of basis vectors (and with respect to the standard basis in $\mC^n$)
can now be obtained by inspection:
$$
A = \begin{bmatrix} 0 & 0 & 0 \\ 1 & 0 & 0 \\ 0 & 0 & 0 \end{bmatrix}, \qquad
B = \begin{bmatrix} 1 & 0 & 0 \\ 0 & 0 & 0 \\ 0 & 1 & 0 \end{bmatrix}, \qquad
C = \begin{bmatrix} 0 & 1 & -1 \\ 0 & 0 & 1 \\ -1 & 0 & 0 \end{bmatrix}.
$$
One can verify directly that indeed $C(\lambda I - A)^{-1}B \doteq T^{-1}$.
\end{example}

Gohberg and Sigal \cite{GS} give a principal part formula in which both the right and the left canonical matrix are supposed to be
chosen in a particular way. In the version of Keldysh \cite{Keldysh1951,Keldysh1971}, the right canonical matrix can be
arbitrarily chosen, but the left canonical matrix needs to be selected in a special way to match the choice of the
right canonical matrix. One may wonder whether it is possible to give a principal part formula in which both 
the left and the right canonical matrix can be arbitrarily chosen. This question is answered in the positive 
below. First, we show an alternative form of the biorthogonality condition.

\begin{proposition} \label{vty}
Let $V \in \cH^{r \times n}$ and $Y \in \cH^{n \times r}$ be left and right canonical matrices, respectively, for $T$ 
at $\lambda_0$. Then the matrix $VTY$ is invertible, and the biorthogonality relation (\ref{biorth}) holds if and only if
\begin{equation} \label{bio}
(VTY)^{-1} \doteq \Delta^{-1}.
\end{equation}
\end{proposition}
\begin{proof}
It has been shown in
Cor.\,\ref{VTYD} that $\Delta^{-1}VTY$ is unimodular;  hence, $VTY$ is invertible. Suppose now 
that (\ref{biorth}) holds, so that there exists a matrix $H \in \cH_0^{r \times r}$ such that $\Delta^{-1}VTY\!\Delta^{-1} 
= \Delta^{-1} + H$. From this one finds, upon multiplying from the left by $(\Delta^{-1}VTY)^{-1}$, that
$$
\Delta^{-1} = (\Delta^{-1}VTY)^{-1}\Delta^{-1} + (\Delta^{-1}VTY)^{-1}H = (VTY)^{-1} + (\Delta^{-1}VTY)^{-1}H
$$
which implies that $\Delta^{-1} \doteq (VTY)^{-1}$. 

Conversely, suppose there exists a matrix $H \in \cH_0^{r \times r}$ such
that $(VTY)^{-1} = \Delta^{-1}+H$. Multiplying from the right by $\Delta$, one finds $I+H\Delta = (VTY)^{-1}\Delta
= (\Delta^{-1}VTY)^{-1}$, so that $\Delta^{-1}VTY = (I+H\Delta)^{-1}$. We can then write
$$
\Delta^{-1} VTY \Delta^{-1} = (I+H\Delta)^{-1} \Delta^{-1} = (I+H\Delta)^{-1}\big( (I+H\Delta)\Delta^{-1} - H \big)
= \Delta^{-1} - (I+H\Delta)^{-1}H 
$$
which shows that $\Delta^{-1} VTY \Delta^{-1} \doteq \Delta^{-1}$.
\end{proof}
\noindent
The following two lemmas will be used in the proof of the main novel result of this paper.

\begin{lemma} \label{onesided}
If $M \in \cH_0^{r \times n}$ is right unimodular, then there exists a unimodular matrix $U_R = [U_R^1 \;\; U_R^2]$
in $\cH_0^{n \times (r+(n-r))}$ such that $M[U_R^1 \;\; U_R^2] = [I_r \;\; 0]$.
\end{lemma}
\begin{proof}
Since $M$ is right unimodular, the matrix $M(\lambda_0) \in \mC^{r \times n}$ has full row rank. Consequently, we 
can choose a matrix $\tM_0 \in \mC^{(n-r)\times n}$ such that the matrix $\hM_0 \in \mC^{n\times n}$ obtained by 
stacking $M(\lambda_0)$ and $\tM_0$ is nonsingular. The matrix $\hM \in \cH_0^{n \times n}$ obtained by
stacking $M$ and $\tM_0$ is then unimodular. Define $U_R = \hM^{-1} \in \cH_0^{n\times n}$, and partition this matrix
as $U_R = [ U_R^1 \;\; U_R^2 ]$ with $U_R^1 \in \cH_0^{n \times r}$. We have
$$
\begin{bmatrix} M \\ \tM_0 \end{bmatrix} \begin{bmatrix} U_R^1 \;\; U_R^2 \end{bmatrix} =
\begin{bmatrix} I_r & 0 \\ 0 & I_{n-r} \end{bmatrix}
$$
so that in particular $M[ U_R^1 \;\; U_R^2 ] = [I_r \;\; 0]$. 
\end{proof}

\begin{lemma} \label{nonsing}
Let $M \in \mC^{n \times n}$, and write $\dim \ker M = r$. Let $N_L \in \mC^{(n-r) \times n}$ and 
$N_R \in \mC^{n \times (n-r)}$ be matrices such that $N_R$ has full column rank,
$\ker N_L \cap \im M = \{0\}$, and $\im N_R \cap \ker M = \{0\}$. Then the matrix 
$N_LMN_R \in \mC^{(n-r)\times(n-r)}$ is nonsingular.
\end{lemma}

\begin{proof}
Suppose $x \in \mC^n$ is such that $N_LMN_R x = 0$. Then $MN_R x \in \im M \cap \ker N_L = \{0\}$, so that 
$MN_Rx=0$. Consequently, $N_Rx \in \ker M \cap \im N_R = \{0\}$, so that $N_Rx = 0$ and hence
$x=0$. 
\end{proof}

\begin{theorem} \label{mainthm}
Let $V \in \cH^{r \times n}$ and $Y \in \cH^{n \times r}$ be left and right canonical matrices, 
respectively, for $T$ at $\lambda_0$. 
\nw{Then $VTY$ is invertible (Cor.\,\ref{VTYD}), and the following statements hold:}
\begin{itemize}
\item[\nw{(i)}] \nw{$\hY := Y(VTY)^{-1}\Delta$ is a right canonical matrix for $T$ at $\lambda_0$ and satisfies 
the biorthogonality condition (\ref{biorth}) with respect to $V$;}
\item[\nw{(ii)}] \nw{$\hV := \Delta(VTY)^{-1}V$ is a left canonical matrix for $T$ at $\lambda_0$ and satisfies 
the biorthogonality condition (\ref{biorth}) with respect to $Y$;}
\item[\nw{(iii)}]  the principal part of $T^{-1}$ at $\lambda_0$ is given by
\begin{equation} \label{main}
T^{-1} \doteq Y(VTY)^{-1}V.
\end{equation}
\end{itemize}
\end{theorem}

\begin{proof}
Since $\Delta^{-1}VT \in \cH_0^{r\times n}$ is right unimodular, it follows from Lemma \ref{onesided} that
there exists a left unimodular matrix $\tY \in \cH_0^{n \times (n-r)}$ such that $\Delta^{-1}VT\tY = 0$.
In particular, $(\Delta^{-1}VT)(\lambda_0) \tY(\lambda_0) = 0$. Because $\tY(\lambda_0)$ has full column rank,
it follows that $\im \tY(\lambda_0) = \ker\, (\Delta^{-1}VT)(\lambda_0)$. Then Prop.\,\ref{crux2} shows
that $\im \tY(\lambda_0) \cap \ker T(\lambda_0) = \{0\}$. Since $Y(\lambda_0)$ is a basis matrix for
$\ker T(\lambda_0)$, this means that the
composite matrix $[Y(\lambda_0) \;\; \tY(\lambda_0)] \in \mC^{n \times n}$ is nonsingular, and hence 
the matrix $[ Y \;\; \tY ] \in \cH_0^{n \times n}$ is unimodular. From the fact that $\Delta^{-1}VT\tY = 0$, it also follows 
that $VT\tY = 0$. Similarly, one can find a right unimodular matrix $\tV \in \cH_0^{(n-r) \times n}$ such that 
$\im T(\lambda_0) \cap \ker \tV(\lambda_0) = \{0\}$, the composite
matrix in $\cH_0^{n \times n}$ formed from $V$ and $\tV$ is unimodular, and $\tV TY = 0$. 
It follows
from Lemma \ref{nonsing} that the matrix $\tV(\lambda_0)T(\lambda_0)\tY(\lambda_0) \in \mC^{(n-r) \times
(n-r)}$ is nonsingular, so that $\tV T\tY \in \cH_0^{(n-r) \times (n-r)}$ is unimodular. From
\begin{equation} \label{uni}
\begin{bmatrix} V \\ \tV \end{bmatrix} T \begin{bmatrix} Y & \tY \end{bmatrix} =
\begin{bmatrix} VTY & 0 \\ 0 & \tV T\tY \end{bmatrix}
\end{equation}
it follows that
$$
\begin{bmatrix}
TY(VTY)^{-1} & T\tY(\tV T \tY)^{-1} 
\end{bmatrix}
= T \begin{bmatrix} Y & \tY \end{bmatrix} \begin{bmatrix}
VTY & 0 \\ 0 & \tV T \tY \end{bmatrix}^{-1} = \begin{bmatrix} V \\ \tV \end{bmatrix}^{-1} \in \cH_0^{n \times n}
$$
which shows that $TY(VTY)^{-1} \doteq 0$. Consequently, the matrix $\hY$ defined in claim (i) satisfies
$T\hY \Delta^{-1} = TY(VTY)^{-1} \doteq 0$. Since $\Delta^{-1}VTY$ is unimodular by Cor.\,\ref{VTYD}, we have
$\hY = Y(\Delta^{-1}VTY)^{-1} \in \cH_0^{n \times r}$, and
$\im \hY(\lambda_0) = \im Y(\lambda_0)$. Hence, it follows from Prop.\,\ref{canmat}
that $\hY$ is a right canonical matrix for $T$ at $\lambda_0$. It is trivial to verify that $\hY$ satisfies the
biorthogonality condition with respect to $V$; indeed, we have in fact $\Delta^{-1}VT\hY \Delta^{-1} = \Delta^{-1}$.
This proves claim (i). Claim (ii) is proved analogously. To prove claim (iii), it now suffices to use claim (ii) together 
with Thm.\,\ref{km}.
\end{proof}
\begin{remark}
The validity of claim (iii) can also be concluded without making use of Thm.\,\ref{km}. Indeed, from (\ref{uni}) we have
\old{
$$
T^{-1} = \begin{bmatrix} Y & \tY \end{bmatrix}
\begin{bmatrix} (VTY)^{-1} & 0 \\ 0 & (\tV T\tY)^{-1} \end{bmatrix}
\begin{bmatrix} V \\ \tV \end{bmatrix} 
= Y(VTY)^{-1}V  + \tY(\tV T\tY)^{-1} \tV.
$$
This proves the claim, since $\tY(\tV T\tY)^{-1} \tV$ is holomorphic in a neighborhood of $\lambda_0$.}
\end{remark}

\begin{example}
We can apply formula (\ref{main}) in the setting of Example \ref{xmpl} with $Y$ and $V$ given by (\ref{Ydef}) and
(\ref{Vdef}) respectively. We have
$$
(VTY)(\lambda) = \begin{bmatrix} \lambda^2 & \lambda^2 \\ 0 & \lambda \end{bmatrix}
$$
so that
$$
(VTY)^{-1}(\lambda) = \begin{bmatrix} \lambda^{-2} & -\lambda^{-1} \\ 0 & \lambda^{-1} \end{bmatrix}.
$$
Using this, one finds
$$
(Y(VTY)^{-1}V)(\lambda) = \begin{bmatrix} \lambda^{-2} & - \lambda^{-1} & 0 \\ 0 & \lambda^{-1} & 0 \\ -\lambda^{-1} & 1 & 0 
\end{bmatrix} \doteq \begin{bmatrix} \lambda^{-2} & - \lambda^{-1} & 0 \\ 0 & \lambda^{-1} & 0 \\ -\lambda^{-1} & 0 & 0 
\end{bmatrix}
$$
in agreement with (\ref{pp}).
\end{example}
\noindent
Item (i) of Thm.\,\ref{mainthm} provides a solution to the construction problem that is implicit in Thm.\,\ref{km}. 

\begin{example}
Consider again Example \ref{xmpl} with $Y$ and $V$ given by (\ref{Ydef}) and (\ref{Vdef}) respectively. We find
$$
\hY(\lambda) = Y(\lambda)(VTY)^{-1}(\lambda) \Delta(\lambda) = \begin{bmatrix} 1 & -1 \\ 0 & 1 \\ -\lambda & \lambda \end{bmatrix}.
$$
This agrees with the solution $\bY$ found in Example \ref{xmpl3}, up to a matrix of the form $H\Delta$ with 
$H \in \cH_0^{3 \times 2}$.
\end{example}

Calculation directly from the formula (\ref{main}) requires inversion of a matrix in $\cH^{r \times r}$, whereas
the method suggested by Thm.\,\ref{km} can be carried out on the basis of solving a linear system over $\mC$,
as illustrated in Example \ref{xmpl3}. To obtain a computation over $\mC$ from (\ref{main}), one can rewrite
the formula as 
\begin{equation} \label{alt}
T^{-1} \doteq Y\Delta^{-1}(VTY\Delta^{-1})^{-1}V
\end{equation} 
and use the fact that, by Cor.\,\ref{VTYD}, the matrix $VTY\Delta^{-1}$
is unimodular. Generally speaking, the coefficients of the power series expansion around $\lambda_0$ of the inverse of 
a unimodular matrix $M = \sum_{j=0}^\infty M_j \chi_0^j$ can be computed recursively from
\begin{equation} \label{recur}
(M^{-1})_0 = M_0^{-1}, \qquad (M^{-1})_k = -M_0^{-1} \sum_{j=1}^k M_j (M^{-1})_{k-j} \quad (k=1,2,,\dots).
\end{equation}
For the purposes of finding the principal part of $T^{-1}$ from (\ref{alt}), it is sufficient to have the power series
coefficients of $(VTY\Delta^{-1})^{-1}$ available for indices up to the order of the pole of $T$ at $\lambda_0$ minus one.
Instead of using (\ref{alt}), one can, of course, also use the mirrored version $T^{-1} = Y(\Delta^{-1}VTY)^{-1}\Delta^{-1}V$;
this may lead to an easier computation.

\begin{example}
Consider once more Example \ref{xmpl} with $Y$ and $V$ given by (\ref{Ydef}) and (\ref{Vdef}) respectively. 
We have
$$
(VTY\Delta^{-1})(\lambda) = \begin{bmatrix} \lambda^2 & \lambda^2 \\ 0 & \lambda \end{bmatrix}
\begin{bmatrix} \lambda^{-2} & 0 \\ 0 & \lambda^{-1} \end{bmatrix} = 
\begin{bmatrix} 1 & \lambda \\ 0 & 1 \end{bmatrix}.
$$
To apply (\ref{alt}), it is sufficient to have the first two coefficients of the power series expansion of the inverse 
of this matrix around $\lambda = 0$ available; these can be computed using (\ref{recur}). One finds
$$
T^{-1} \doteq \begin{bmatrix} 1 & 0 \\ 0 & 1 \\ -\lambda & 0 \end{bmatrix} 
\begin{bmatrix} \lambda^{-2} & 0 \\ 0 & \lambda^{-1} \end{bmatrix} 
\begin{bmatrix} 1 & -\lambda \\ 0 & 1 \end{bmatrix} 
\begin{bmatrix} 1 & 0 & 0 \\ 0 & 1 & 0 \end{bmatrix} = 
\begin{bmatrix} \lambda^{-2} & - \lambda^{-1} & 0 \\ 0 & \lambda^{-1} & 0 \\ -\lambda^{-1} & 1 & 0 
\end{bmatrix}
$$
in agreement with (\ref{pp}). Using the mirrored version of (\ref{alt}) makes for a slightly simpler
computation, since $\Delta^{-1}VTY$ is in fact constant in the case at hand.
\end{example}

\nw{Thm.\,\ref{mainthm}} shows that both the condition $VTY = \Delta$ (as in Thm.\,\ref{gs}, see Remark \ref{VTYw})
and the condition $(VTY)^{-1} \doteq \Delta^{-1}$ (as in Thm.\,\ref{keldysh} and Thm.\,\ref{km}, see Prop.\,\ref{vty})
are sufficient for Keldysh's principal part formula (\ref{rpp}) to hold. 

\section{The semisimple case}

The root of $T$ at $\lambda_0$ is said to be \emph{semisimple} when there are no partial multiplicities larger than 1,
i.e., $m_1 = 1$. As noted below Thm.\,\ref{Smith}, this condition is necessary and sufficient 
for $T^{-1}$ to have a simple pole (i.e., a pole of order 1) at $\lambda_0$. In several applications, it
is of interest to characterize semisimplicity in terms of root spaces associated to $T$ at $\lambda_0$, and to have
an expression for the residue of $T^{-1}$ in the semisimple case. Conditions for semisimplicity are not discussed
explicitly by Keldysh \cite{Keldysh1951,Keldysh1971}, but he does prove \cite[p.\,21]{Keldysh1971} that the pole order
is equal to the maximal multiplicity of root functions. 
A condition for semisimplicity in terms of the first two coefficients of the power series expansion of $T$ around $\lambda_0$
can be derived as follows.

\begin{proposition} \label{char1}
The root of $T$ at $\lambda_0$ is semisimple if and only if\,\footnote{\nw{Given a mapping $A\cn X \rightarrow Y$ and a subset
$Y_0$ of $Y$, the notation $A^{-1}(Y_0)$ indicates the \emph{inverse image} of $Y_0$ under the mapping $A$, i.e., $A^{-1}(Y_0) :=
\{ x \in X \mid Ax \in Y_0 \}$. Invertibility of $A$ is not assumed.}}
\begin{equation} \label{cond1}
\ker T(\lambda_0) \cap (T'(\lambda_0))^{-1}(\im T(\lambda_0)) = \{0\}.
\end{equation}
\end{proposition}

\begin{proof}
By Prop.\,\ref{jl}, we have
$m_1=1$ if and only if $\ell_2=0$, where $\ell_2 = \dim L_2$, and $L_2$ is defined by (\ref{LjT}).
According to Prop.\,\ref{Ljchar}, the subspace $L_2$ consists of all vectors $y^0$ that satisfy (i)
$T(\lambda_0)y^0 = 0$, and (ii) there exists $y^1 \in \mC^n$ such that $T(\lambda_0)y^1 + T'(\lambda_0)y^0 = 0$. 
In other words, $L_2$ is given by the expression at the left hand side of (\ref{cond1}). 
\end{proof}
\noindent
Alternative characterizations follow from this via the following lemma.

\begin{lemma}
Let matrices $M, N \in \mC^{n\times n}$ be given. Write $r = \dim \ker M$, and let 
$Z_L \in \mC^{r \times n}$ and 
$Z_R \in \mC^{n \times r}$ be basis matrices for the left and right null spaces of $M$ respectively. The following 
statements are equivalent.
\begin{itemize}
\item[(i)] $\ker M \cap N^{-1}(\im M) = \{0\}$.
\item[(ii)] $\ker M \cap \ker N = \{0\}$ and $\im M \cap N(\ker M) = \{0\}$.
\item[(iii)] $\mC^n = \im M \oplus N(\ker M)$.
\item[(iv)] The matrix $Z_L N Z_R \in \mC^{r \times r}$ is nonsingular.
\end{itemize} 
\end{lemma}
\begin{proof}
Assume that (i) holds. Then $\ker N \cap \ker M = \{0\}$, since $\ker N \subset N^{-1}(\im M)$. Take 
$x \in \im M \cap N(\ker M)$, and let $y \in \ker M$ be such that $x=Ny$. We then have $y \in 
\ker M \cap N^{-1}(\im M)$, so that $y=0$ and hence $x=0$. This proves that (i) implies (ii). Now, assume that
(ii) holds. Take $x \in \ker M \cap N^{-1}(\im M)$. Then $Nx \in \im M \cap N(\ker M)$ so that $Nx = 0$.
It follows that $x \in \ker N \cap \ker M$; hence, $x=0$. This shows that (ii) implies (i). Equivalence of (ii) and (iii)
follows by noting that
\begin{align*}
\dim(\im M + N(\ker M)) & = \dim \im M + \dim N(\ker M)) - \dim(\im M \cap N(\ker M)) \\
& = \dim \im M + \dim \ker M - \dim(\ker N \cap \ker M) - \dim(\im M \cap N(\ker M)) \\
& = n - \big(\!\dim(\ker N \cap \ker M) + \dim(\im M \cap N(\ker M))\big).
\end{align*}
Finally, observe that
\begin{align*}
\dim \im Z_L NZ_R & = \dim N(\im Z_R)  - \dim(\ker Z_L \cap N(\im Z_R)) \\
& = \dim N(\ker M) - \dim(\im M \cap N(\ker M)) \\
& = r - \big(\!\dim(\ker N \cap \ker M) + \dim(\im M \cap N(\ker M))\big).
\end{align*} 
Since the matrix $Z_L NZ_R$ has size $r \times r$, it follows that (ii) and (iv) are equivalent.
\end{proof}
\noindent

\begin{proposition} \label{semcon}
Let $Z_L$ and $Z_R$ be basis matrices for the left and right null spaces of $T(\lambda_0)$ respectively.
Each of the following conditions is necessary and sufficient for semisimplicity of the root of $T $ at $\lambda_0$.
\begin{itemize}
\item[(i)] $\mC^n = \im T(\lambda_0) \oplus T'(\lambda_0)(\ker T(\lambda_0))$.
\item[(ii)] $\ker T(\lambda_0) \cap \ker T'(\lambda_0) = \{0\}$ and 
$\im T(\lambda_0) \cap T'(\lambda_0)(\ker T(\lambda_0)) = \{0\}$.
\item[(iii)] The matrix $Z_L T'(\lambda_0) Z_R \in \mC^{r \times r}$ is nonsingular.
\end{itemize}
\end{proposition}

The necessity of condition (iii) can also be derived from Cor.\,\ref{VTYD}. Indeed, in the semisimple case, the 
corollary implies that the matrix $\chi_0^{-1}Z_RTZ_L \in \cH_0^{r \times r}$ is unimodular. Since 
$(\chi_0^{-1}Z_RTZ_L)(\lambda_0) = Z_RT'(\lambda_0)Z_L$, it follows that $Z_RT'(\lambda_0)Z_L$ is nonsingular. 

If semisimplicity holds, the principal part of $T^{-1}$ at $\lambda_0$ is completely specified by the residue, denoted by 
$\text{\rm Res}(T^{-1};\lambda_0)$. A formula for the residue can be obtained from Keldysh's theorem by
specializing the notions of left and right canonical matrices to the semisimple case. At a semisimple root,
canonical matrices can be constructed using constant vectors only. A matrix in $\mC^{r \times n}$
($\mC^{n\times r}$) is left (right) canonical if and only if it is a basis matrix for the left (right) null space of
$T(\lambda_0)$. 
From Thm.\,\ref{keldysh}, one then has the following corollary.

\begin{corollary} \label{kelcor}
If the root of $T$ at $\lambda_0$ is semisimple, then, given any basis matrix $Z_R \in \mC^{n\times r}$ for the
right null space of $T(\lambda_0)$, there exists a uniquely determined basis matrix $Z_L \in \mC^{r \times n}$ 
for the left null space of $T(\lambda_0)$ such that the residue 
of $T^{-1}$ at $\lambda_0$ is given by $Z_R Z_L$. The matrix $Z_L$ satisfies
\begin{equation} \label{bicst}
Z_L T'(\lambda_0)Z_R = I.
\end{equation}
\end{corollary}
\begin{proof}
Let $Z_R$ be a basis matrix for $\ker T(\lambda_0)$. Thm.\,\ref{keldysh} implies that there exists a uniquely
determined constant matrix $Z_L$
such that $T^{-1} \doteq Z_R\Delta^{-1}Z_L = \chi_0^{-1} Z_R Z_L$. This shows that the residue of $T^{-1}$ at 
$\lambda_0$ is given by $Z_R Z_L$. Since $Z_L T(\lambda_0) Z_R=0$, it follows from the biorthogonality 
condition (\ref{biorth}) that 
$$\
\chi_0^{-1}I \doteq \chi_0^{-2}Z_L T Z_R \doteq \chi_0^{-1}Z_L T'(\lambda_0) Z_R.
$$
This implies (\ref{bicst}).
\end{proof}
One can now reason as follows. Let a basis matrix $Z_R \in \mC^{n \times r}$ for $\ker T(\lambda_0)$ be given, 
and let $Z_L \in \mC^{r \times n}$ be an \emph{arbitrary} 
left basis matrix for the left null space of $T(\lambda_0)$. Let $\tZ_L$ denote the left basis matrix obtained from 
Cor.\,\ref{kelcor}. Since $Z_L$ and $\tZ_L$ are both left basis matrices, there exists a nonsingular matrix $M \in
\mC^{r \times r}$ such that $\tZ_L = MZ_L$. Equation (\ref{bicst}) then implies that $MZ_L T'(\lambda_0)Z_R = I$, 
and hence $M = (Z_L T'(\lambda_0)Z_R)^{-1}$.
It follows that $\tZ_L = (Z_L T'(\lambda_0) Z_R)^{-1}Z_L$. This leads to the next corollary.

\begin{corollary} \label{kelcor2}
Suppose that the root of $T$ at $\lambda_0$ is semisimple. Let $Z_L \in \mC^{r \times n}$ and 
$Z_R \in \mC^{n \times r}$ be basis matrices for the left and the right null spaces of $T(\lambda_0)$ respectively. 
Then the matrix $Z_L T'(\lambda_0) Z_R \in \mC^{r \times r}$ is nonsingular, and the residue of $T^{-1}$ at 
$\lambda_0$ is given by 
\begin{equation} \label{resform}
\text{\rm Res}(T^{-1};\lambda_0) = Z_R(Z_LT'(\lambda_0)Z_R)^{-1}Z_L.
\end{equation}
\end{corollary}
\vskip2mm\noindent

Semisimplicity conditions and corresponding residue formulas have a somewhat complicated history.
As noted above, conditions for semisimplicity are not given explicitly in \cite{Keldysh1951,Keldysh1971}.
For polynomial matrices $T$ with nonsingular leading coefficient matrix, Lancaster \cite[Thm.\,4.3]{Lancaster}
proves that a necessary condition for the pole of $T^{-1}$ at $\lambda_0$
to be of first order is that there exist basis matrices $Z_L$ and $Z_R$ for the left and right null spaces of $Z(\lambda_0)$
such that $Z_L T'(\lambda_0) Z_R = I$, and that the residue is then given by $Z_R Z_L$. The proof uses 
a linearization technique.\footnote{Linearization is meant here in the sense that a linear pencil of matrices is used
to represent the given polynomial matrix $T$. The assumption that the highest-order coefficient matrix of $T$ is
nonsingular is expedient in this context since it allows a straightforward generalization of the companion form,
a standard tool in the theory of higher-order scalar differential equations.} 
Working in an infinite-dimensional setting and also using a linearization technique, Howland \cite{Howland} shows that (i) and 
(ii) in Prop.\,\ref{semcon} are both necessary and sufficient for semisimplicity. 
This extends earlier work of Steinberg \cite{Steinberg}. 
Neither Howland nor Steinberg gives a formula for the residue, and both do not refer to
Keldysh's work. \nw{For the matrix case, Howlett \cite{Howlett} proves sufficiency of another condition equivalent
to (i)--(iii) of Prop.\,\ref{semcon}, again without residue formula or mention of Keldysh. Also i}n the matrix setting, 
Schumacher \cite{ESTEC,MTNS} gives condition (iii) as well as the residue formula (\ref{resform}), with a derivation
on the basis of the Smith form. The 1951 paper by Keldysh is referred to, unfortunately with the
incorrect claim that the normalization factor $(Z_L T'(\lambda_0) Z_R)^{-1}$ is missing in that paper.\footnote{%
This calls for an explanation. When writing the cited publications, I was unaware of \cite{Keldysh1971} and had 
no direct access to \cite{Keldysh1951}, but I did know about the existence of the latter paper from \cite{GS}.
However, Gohberg and Sigal cite Keldysh's theorem in a version in which the biorthogonality condition is omitted.
From the fact that the normalization factor does not appear in their paper, I wrongly concluded that it had 
also not been given by Keldysh. While the residue formula (\ref{resform}) is not stated explicitly in \cite{Keldysh1951}, 
using the biorthogonality condition one can, as shown above, derive it from Keldysh's theorem, including the 
normalization factor.}
Condition (iii) also appears in \cite[Prop.\,1.7.3]{MM},\footnote{The proposition in fact states as a condition that 
there should exist basis matrices $Z_L$ and $Z_R$ such that $Z_L T'(\lambda_0) Z_R$ is nonsingular. Since left (right)
basis matrices are defined uniquely up to left (right) multiplication by a nonsingular matrix, it is 
clear that the condition holds for a given pair of basis matrices $(Z_L,Z_R)$ if and only if it holds for all such pairs.}
with a derivation on the basis of Keldysh's theorem. \nw{A proof of sufficiency of (iii), stated for a specific application 
in vector time series analysis and formulated in the corresponding language, is given by Johansen \cite{Johansen} who also similarly arrives 
at the residue formula (\ref{resform}).}

The root of $T$ at $\lambda_0$ is said to be \emph{simple} when in addition to $m_1=1$ we also have $r=1$, i.e.,
the algebraic multiplicity of the pole is equal to 1. In this case, Keldysh's theorem Thm.\,\ref{keldysh} implies that,
for any right null vector $y^0$ of $T(\lambda_0)$, there exists a left null vector $v^0$ of $T(\lambda_0)$ such that
$T^{-1} \doteq \chi_0^{-1} y^0 v^0$. Since left null vectors are determined up to a nonzero scalar factor, it is
immediate from the biorthogonality condition $\chi_0^{-2} v^0Ty^0 \doteq \chi_0^{-1}$ that the vector $v^0$ 
required in Keldysh's theorem should be chosen such that $v^0T'(\lambda_0)y^0 = 1$. Beyn \cite{Beyn} gives this result 
and refers to it as Keldysh's theorem for simple eigenvalues.

The residue formula (\ref{resform}) for first-order poles is akin to the principal part formula (\ref{main}). One can
derive the former from the latter in the following way. Let it be given that the root of $T$ at $\lambda_0$ is 
semisimple, and let $Z_L$ and $Z_R$ be basis matrices for the left and right null spaces of $T(\lambda_0)$ respectively. 
We have
$$
(Z_RTZ_L)^{-1} = (\chi_0^{-1}Z_RTZ_L)^{-1}\chi_0^{-1} \doteq (Z_RT'(\lambda_0)Z_L)^{-1} \chi_0^{-1}.
$$
Applying (\ref{main}), one arrives at the residue formula (\ref{resform}).

As a final comment, we note that, while the semisimplicity condition in Prop.\,\ref{semcon}(iii) and the residue formula 
(\ref{resform}) can be obtained from Keldysh's theorem as shown above, a short derivation from the Smith
form is possible. The proof below summarizes the argument in \cite{ESTEC,MTNS}.

\begin{proposition} \label{residuethm}
Let $Z_L$ and $Z_R$ be basis matrices for the left and right null spaces of $T(\lambda_0)$ respectively.
The root of $T$ at $\lambda_0$ is semisimple if and only if the matrix $Z_L T'(\lambda_0) Z_R \in \mC^{r \times r}$ is 
nonsingular, and in that case the residue of $T^{-1}$ at $\lambda_0$ is given by (\ref{resform}).
\end{proposition}

\begin{proof}
Let $\tT = U_L T U_R$ where $U_L,\, U_R \in \cH_0^{n\times n}$ are unimodular.
Since $\tT(\lambda_0) = U_L(\lambda_0)T(\lambda_0)U_R(\lambda_0)$ where $U_L(\lambda_0)$ and 
$U_R(\lambda_0)$ are nonsingular, we have $\dim \ker \tT(\lambda_0) = \dim \ker T(\lambda_0) = r$. The pole 
of $\tT^{-1}$ at $\lambda_0$ is of first order if and only if the same holds for $T^{-1}$, and in that case we have
\begin{equation} \label{restrafo}
\text{\rm Res}(\tT^{-1};\lambda_0) = (U_R(\lambda_0))^{-1} \text{\rm Res}(T^{-1};\lambda_0) (U_L(\lambda_0))^{-1}.
\end{equation}
Let $Z_L\, (\tZ_L) \in \mC^{r\times n}$ and
$Z_R\, (\tZ_R) \in \mC^{n\times r}$ be basis matrices for the left and right null spaces
of $T(\lambda_0)$ ($\tT(\lambda_0)$). Since $U_R(\lambda_0)\tZ_R$ 
and $Z_R$ are both basis matrices for $\ker T(\lambda_0)$,
there exists a nonsingular matrix $M_R \in \mC^{r \times r}$ such that 
$U_R(\lambda_0) \tZ_R = Z_R M_R$. Likewise, there exists a nonsingular matrix $M_L \in \mC^{r \times r}$ such that 
$\tZ_L U_L(\lambda_0) = M_L Z_L$. Hence,
\begin{align*}
\tZ_L \tT'(\lambda_0) \tZ_R  & = \tZ_L \big( U_L'(\lambda_0)T(\lambda_0)U_R(\lambda_0) + U_L(\lambda_0)T'(\lambda_0)U_R(\lambda_0)
+ U_L(\lambda_0)T(\lambda_0)U_R'(\lambda_0) \big) \tZ_R \\
& = M_L Z_L T'(\lambda_0) Z_R M_R.
\end{align*}
It follows that the matrix $\tZ_L \tT'(\lambda_0) \tZ_R \in \mC^{r \times r}$ is nonsingular if and only if the same 
is true for $Z_L T'(\lambda_0) Z_R$. If this holds, then it follows by comparing (\ref{restrafo}) with
\begin{equation} \label{restr}
\tZ_R(\tZ_L \tT'(\lambda_0) \tZ_R)^{-1}\tZ_L = 
(U_R(\lambda_0))^{-1} Z_R(Z_L T'(\lambda_0) Z_R )^{-1}Z_L(U_L(\lambda_0))^{-1}
\end{equation}
that the residue formula (\ref{resform}) is valid for 
$T$ if and only if it is valid for $\tT$. To prove the claims of the proposition, it therefore suffices to verify
them for a matrix in Smith form; basis matrices for left and right null spaces may be arbitrarily chosen. 
Therefore, suppose the matrix $T$ is given by
$$
T(\lambda) = \begin{bmatrix} \Delta(\lambda) & 0 \\ 0 & I_{n-r} \end{bmatrix} \qquad\qquad (\lambda \in \Omega)
$$
with $\Delta = \diag(\chi_0^{m_1}, \dots, \chi_0^{m_r})$ as in (\ref{deltadef}). Basis matrices for the left and right 
null spaces of $T$ at $\lambda_0$ are given by 
$$
Z_L = \begin{bmatrix} I_r & 0 \end{bmatrix}, \qquad Z_R = \begin{bmatrix} I_r \\ 0 \end{bmatrix}.
$$
We have $Z_L T'(\lambda_0) Z_R = \Delta'(\lambda_0)$. This matrix is nonsingular if and only if $m_i=1$
for all $i=1,\dots,r$, i.e., if and only if the root of $T$ at $\lambda_0$ is semisimple. If this holds, then
$\Delta = \chi_0 I_r$ so that $Z_L T'(\lambda_0) Z_R = Z_L Z_R = I_r$, and
$$
T^{-1} = \begin{bmatrix} \chi_0^{-1} I_r & 0 \\ 0 & I_{n-r} \end{bmatrix} = 
\begin{bmatrix} 0 & 0 \\ 0 & I_{n-r} \end{bmatrix} + \chi_0^{-1} 
\begin{bmatrix}  I_r & 0 \\ 0 & 0 \end{bmatrix}
\doteq \chi_0^{-1} Z_R(Z_LT'(\lambda_0)Z_R)^{-1}Z_L.
$$
\vskip-4mm
\end{proof}

\section*{Acknowledgement}
I would like to thank Brendan Beare and two anonymous reviewers for their comments on an earlier version of this paper.

\begin{singlespace}
\bibliographystyle{plainnat}
\bibliography{semisimple}
\end{singlespace}
\end{document}